\documentclass[smallextended]{svjour3}       
\smartqed  

\usepackage{amsfonts,amssymb,amsmath,amsopn,float}
\usepackage{xifthen}
\usepackage{graphicx}
\usepackage{tikz}
\usepackage{pgfplots}
\usepgfplotslibrary{units}
\usepackage{epstopdf}
\usepackage[utf8]{inputenc}
\usepackage{mathtools}
\usepackage{url}
\usepackage{subcaption}

\usepackage{multicol}

\usepackage[pdftex,hypertexnames=true,plainpages=false,
naturalnames,colorlinks=true,linkcolor=blue,bookmarks]{hyperref}

\newcommand{\bbR}{\mathbb{R}}

\newcommand{\Ltwo}[1]{%
	\ifthenelse{\equal{#1}{}}{L^2}{L^2(#1)}%
}

\newcommand{\Ltwoz}[1]{%
	\ifthenelse{\equal{#1}{}}{L^2_0}{L^2_0(#1)}%
}

\newcommand{\Cone}[1]{%
	\ifthenelse{\equal{#1}{}}{C^{1}}{C^{1}(#1)}%
}

\newcommand{\Conez}[1]{%
	\ifthenelse{\equal{#1}{}}{C^{1}_{0}}{C^{1}_{0}(#1)}%
}

\newcommand{\Ctwo}[1]{%
	\ifthenelse{\equal{#1}{}}{C^{2}}{C^2(#1)}%
}

\newcommand{\Ctwoz}[1]{%
	\ifthenelse{\equal{#1}{}}{C^{2}_{0}}{C^{2}_{0}(#1)}%
}

\newcommand{\Cholder}[1]{%
	\ifthenelse{\equal{#1}{}}{C^{0,\gamma}}{C^{0,\gamma}(#1)}%
}

\newcommand{\Cholderz}[1]{%
	\ifthenelse{\equal{#1}{}}{C^{0,\gamma}_{0}}{C^{0,\gamma}_{0}(#1)}%
}

\newcommand{\bolds}[1]{\boldsymbol{#1}}

\newcommand{\bb}{\bolds{b}}
\newcommand{\bc}{\bolds{c}}

\newcommand{\be}{\bolds{e}}

\newcommand{\bn}{\bolds{n}}

\newcommand{\bt}{\bolds{t}}
\newcommand{\bu}{\bolds{u}}
\newcommand{\bv}{\bolds{v}}
\newcommand{\bw}{\bolds{w}}
\newcommand{\bx}{\bolds{x}}
\newcommand{\by}{\bolds{y}}

\newcommand{\sautoref}[2]{\hyperref[#2]{#1 \ref*{#2}}}

\newtheorem{hypothesis}{Hypothesis}

\begin{document}
	
\title{Classic dynamic fracture recovered as the limit of a nonlocal peridynamic model: The single edge notch in tension \thanks{This material is based upon work supported by the U. S. Army Research Laboratory and the U. S. Army Research Office under contract/grant number W911NF1610456.}}

\titlerunning{Classic dynamic fracture recovered as the limit of a nonlocal peridynamic model}  

\author{Robert Lipton
	\and
	Prashant K. Jha}

\institute{Robert Lipton\at
	Department of Mathematics,\\
	and Center for Computation and Technology,\\
	Louisiana State University,\\ 
	Baton Rouge, LA 70803\\
	{Orcid: https://orcid.org/0000-0002-1382-3204}\\
	\email{lipton@lsu.edu}    
	\and
	Prashant K. Jha\at
	Oden Institute for Computational Engineering and Sciences,\\
	The University of Texas at Austin,\\
	Austin, TX 78712\\
	{Orcid: https://orcid.org/0000-0003-2158-364X}\\
	\email{pjha@utexas.edu}    
}

\date{}
\maketitle
\begin{abstract}
A simple nonlocal field theory of peridynamic type is applied to model brittle fracture.
The fracture evolution is shown to converge in the limit of vanishing nonlocality to classic plane elastodynamics with a running crack. The kinetic relation for the crack is recovered directly from the nonlocal model in the limit of vanishing nonlocality. We carry out our analysis for a single crack in a plate subject to mode one loading. The convergence is corroborated by numerical experiments. 
\end{abstract}

  



\section{Introduction}

Fracture can be viewed as a collective interaction across large and small length scales. With the application of enough stress or strain to a brittle material, atomistic scale bonds will break, leading to fracture of the macroscopic specimen.  The appeal of  nonlocal peridynamic models is that fracture appears as an emergent phenomena generated by the underlying field theory eliminating the need for supplemental kinetic relations describing crack growth. The deformation field inside the body for points $\bx$  at time $t$ is written $\bu(\bx,t)$.
The perydynamic model is described simply by the balance of linear momentum of the form
\begin{equation}\label{eqn-uf}
\begin{aligned}
   \rho{{\bu}_{tt}}(\bx,t) = \int_{{\mathcal H}_{\epsilon}({\bx})} {\bolds{f}}(\by,\bx)\;d\by + \bb(\bx,t)
\end{aligned}
\end{equation}
where $\mathcal{H}_{\epsilon({\bx})}$ is a neighborhood of $\bx$,  $\rho$ 
is the density, $\bb$ is the body force 
density field, and $\bolds{f}$ is a material-dependent constitutive law
that represents the force density that a point $\by$ inside the neighborhood exerts on $\bx$ as a result of the deformation field.
The radius $\epsilon$ of the neighborhood is referred to as the \emph{horizon}. 
Here all points satisfy the same basic field equations \eqref{eqn-uf}. This approach to fracture modeling was introduced in \cite{CMPer-Silling} and \cite{States}.  
The displacement fields and fracture evolution predicted by nonlocal models {should agree} with the established theory of dynamic fracture mechanics when the length scale of non-locality is sufficiently small. This phenomena can be seen in simulations, see for example, \cite{ParksTrasketal}, \cite{Bobaru 2015}, and \cite{silling05} .

\begin{figure} 
\centering
\begin{tikzpicture}[xscale=0.60,yscale=0.60]

%
%
\draw [-,thick] (-2,0.05) -- (-1.3,0.05);

\draw [-,thick] (-2,-0.05) -- (-1.3,-0.05);


\draw [-,thick] (-1.3,0.05) to [out=0, in=90] (-1.25 ,0);

\draw [-,thick] (-1.25,0) to [out=-90, in=0] (-1.3 ,-0.05);


%
%
\draw [-,thick] (-2,3) -- (-2,0.05);
\draw [-,thick] (-2,-3) -- (-2,-0.05);
\draw [-,thick] (-2,-3) -- (2,-3);
\draw [-,thick] (2,3) -- (2,-3);
\draw [-,thick] (2,3) -- (-2,3);
%
%

\draw[->,thick] (-2.,3.0) -- (-2,3.80);

\draw[->,thick] (-1,3.0) -- (-1,3.80);

\draw[->,thick] (0.0,3.0) -- (0.0,3.80);

\draw[->,thick] (1,3.0) -- (1,3.80);

\draw[->,thick] (2.0,3.0) -- (2,3.80);


\draw[->,thick] (-2.0,-3.0) -- (-2,-3.80);

\draw[->,thick] (-1,-3.0) -- (-1,-3.80);

\draw[->,thick] (0.0,-3.0) -- (0.0,-3.80);

\draw[->,thick] (1,-3.0) -- (1,-3.80);

\draw[->,thick] (2.0,-3.0) -- (2,-3.80);

\end{tikzpicture} 
\caption{{ \bf Single-edge-notch}}
 \label{singlenotch}
\end{figure}
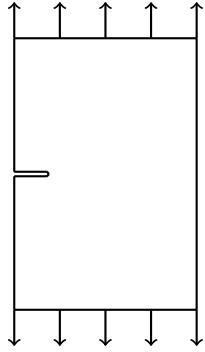

In this paper we theoretically examine the predictions of the nonlocal theory in the limit of vanishing non-locality. We examine a class of peridynamic models with nonlocal forces derived from double well potentials. see \cite{CMPer-Lipton3},  \cite{CMPer-Lipton}.  We theoretically investigate the limit of these evolutions as the length scale $\epsilon$ of nonlocal interaction goes to zero. We are able to describe the interaction between the crack and the surrounding displacement field of intact material in this limit. Here all information on this limit is obtained from what is known from the nonlocal peridynamic model for $\epsilon>0$. We consider a single edge notch specimen as given in Figure \ref{singlenotch}.  For small strains the nonlocal force is linearly elastic but for larger strains the  force begins to soften and then approaches zero after reaching a critical strain. Because of this force vs. strain behavior this type of model is called a cohesive model.

Previous work has addressed the convergence of the cohesive  fracture model to classic local brittle fracture for dynamic free crack propagation with multiple interacting cracks  \cite{CMPer-Lipton3},  \cite{CMPer-Lipton}, \cite{CMPer-JhaLipton}. There it is shown that the nonlocal cohesive  evolution converges to an evolution of sharp cracks with bounded Griffith fracture energy satisfying the linear elastic wave equation off the cracks. However the explicit interaction between the sharp crack and intact material remains to be described in the local limit. In this work we describe in an explicit way the limiting interaction between the sharp crack and surrounding material. Here we  pass to the limit in the nonlocal model to recover the limiting dynamic interaction of the sharp crack with the surrounding intact material. 
A distinguishing feature of the cohesive nonlocal model is that the fracture toughness is the same for all horizons $\epsilon>0$. It is shown here that fracture evolutions are mathematically well posed for every $\epsilon>0$ and that as $\epsilon\rightarrow 0$ the nonlocal evolution converges to the dynamic brittle fracture model given by:
\begin{itemize}
\item
Balance of linear momentum described by the linear elastic wave equation away from the crack.

\item Zero traction on the crack lips.

\item The classic kinetic relation for crack tip velocity implicitly given by equating the dynamic stress intensity factor with the energy dissipation per unit extension of the crack.
\end{itemize}
However in this paper the kinetic relation for crack tip velocity is {\em not derived from the power balance postulate of} \cite{Mott} but  instead is recovered from  the nonlocal model \eqref{energy based model2} directly by taking the $\epsilon=0$ limit in the nonlocal power balance {Proposition} \ref{changinternal} as shown in  Proposition \ref{limitspowerflow}. The kinetic relation derived here follows from an  explicit formula for the time rate of change of internal energy inside a domain containing the crack tip, see Proposition \ref{limitspowerflow}. In this way we recover  the modern dynamic fracture model  developed and described in \cite{Freund}, \cite{RaviChandar}, \cite{Anderson}, \\
\cite{Slepian} but without using the power balance postulate and instead  using the cohesive dynamics based on double well potentials \eqref{energy based model2} directly.  We note further that the limiting classic local fracture problem is hard to simulate directly, this is because the crack velocity at the crack tip is directly coupled to the wave equation off the crack and vice versa. On the other hand this coupling between intact material and crack is handled autonomously in the nonlocal model and numerical simulation is straight forward.  For a-priori convergence rates of finite difference and finite element implementations of the nonlocal model given here  see \cite{CMPer-JhaLipton}, \cite{CMPer-JhaLipton8}  \cite{CMPer-Lipton2}, \cite{CMPer-JhaLipton3}. 

The analysis used in this paper relies in part on the earlier analysis of \cite{CMPer-Lipton} but also requires new compactness methods specifically suited to the balance of momentum for nonlocal - nonlinear operators, see section \ref{s:proofofconvergence}. These methods give the zero traction condition on the crack lips for the fracture model in the local limit. The explicit formula for the time rate of change of the internal energy for a domain containing  the crack tip follows by passing to the $\epsilon\rightarrow 0$ limit in an identity that is obtained using a new type of divergence theorem for nonlocal operators, see section \ref{s:energyrate}. The kinetic relation for the crack tip velocity follows directly from the formula for the time rate of change of internal energy. 

The paper is organized as follows: In section \ref{s:nonlocal dynamics} we describe the nonlocal constitutive law derived from a double well potential and present the nonlocal boundary value problem describing crack evolution. Section \ref{toughelastic} outlines how the fracture toughness and elastic properties of a material are contained in the description of the double well potential. Section \ref{prelim} outlines the preliminary convergence results necessary for the analysis. Section \ref{ss:crackwaveinteraction} provides the principle results of the paper and describes the convergence of the nonlocal crack evolution to the local dynamic fracture evolution described in \cite{Freund}, \cite{RaviChandar}, \cite{Anderson}, \cite{Slepian}. The hypotheses on the emergence and nature of the zone where the force between points decreases with increasing strain follows from the symmetry of the loading and domain and are  corroborated by the numerical simulations in section \ref{numerical}. The existence and uniqueness of the nonlocal evolution is established in section \ref{existenceuniqueness}. 
 The relation between crack set and jump set for the limit evolution is proved in section \ref{lengthjump}. The proof of convergence is given in section \ref{s:proofofconvergence}. The time rate of energy increase inside a region containing the crack tip for the nonlocal model is given in section \ref{s:energyrate} and follows from a new nonlocal divergence theorem given in this section. The kinetic relation for crack tip motion for classic dynamic fracture mechanics is shown to follow directly from the nonlocal model and is derived in section \ref{s:Jintegral}. We summarize results in the conclusion section \ref{s:conclusions}.

\section{Nonlocal Dynamics}\label{s:nonlocal dynamics}

In this section we formulate the nonlocal dynamics as an initial boundary value problem with traction boundary conditions.  We begin by introducing the nonlocal force defined in terms of a double well potential. Here all quantities are non-dimensional.  Define the rectangle $R=\{0<x_1<a;\,-b/2<x_2<b/2\}$ and we will consider the plane strain problem with a thin notch denoted by $C$  of thickness $2d$ and total length $\ell(0)$ with a circular tip. This is described by $\{0\leq x_1\leq\ell(0)-d;\,-d < x_2<d\}\cup\{\ell(0)-d\leq x_1\leq\ell(0)-d+\sqrt{d^2-x_2^2};-d< x_2< d,\}$ originating on the left side of the rectangle. 
The domain for the peridynamic evolution is given by $D=R\setminus C$ and the domain $D$ corresponds to a single edge notch specimen, see Figure \ref{singlenotch}.  In this treatment we will assume small (infinitesimal) deformations so that the displacement field $\bu: D\times [0,T] \to \bbR^2$ is small compared to the size of $D$ and the deformed configuration is the same as the reference configuration.  We have $\bu=\bu(\bx,t)$ as a function of space and time but will  suppress the $\bx$ dependence when convenient and write $\bu(t)$. The tensile strain $S$ between two points $\bx,\by$ in $D$ along the direction $\be_{\by-\bx}$ is defined as
\begin{align}\label{strain}
S(\by,\bx,\bu(t))=\frac{\bu(\by,t)-\bu(\bx,t)}{|\by-\bx|}\cdot \be_{\by-\bx},
\end{align}
where $ \be_{\by-\bx}=\frac{\by-\bx}{|\by-\bx|}$ is a unit vector and ``$\cdot$'' is the dot product. The influence function $J^\epsilon(|\by-\bx|)$ is a measure of the influence that the point $\by$ has on $\bx$. Only points inside the horizon can influence $\bx$ so $J^\epsilon(|\by-\bx|)$ nonzero for $|\by - \bx| < \epsilon$ and zero otherwise. We take $J^\epsilon$ to be of the form: $J^\epsilon(|\by - \bx|) = J(\frac{|\by - \bx|}{\epsilon})$ with $J(r) = 0$ for $r\geq 1$ and $0\leq J(r)\leq M < \infty$ for $r<1$. 

\subsection{The class of nonlocal potentials}
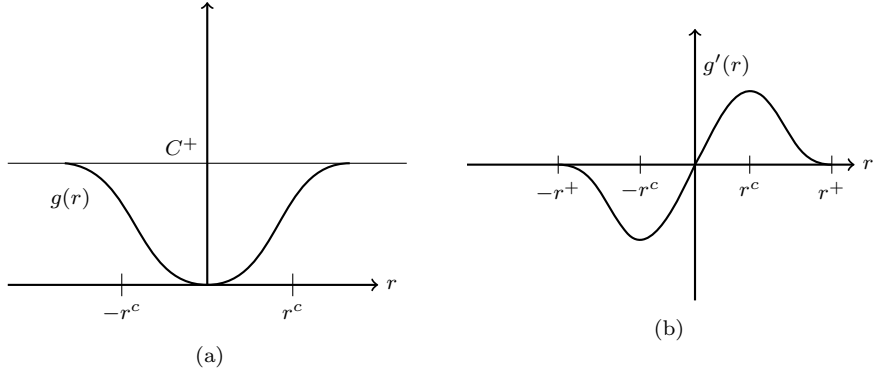
\begin{figure}
    \centering
    \begin{subfigure}{.45\linewidth}
        \begin{tikzpicture}[xscale=0.75,yscale=0.75]
		    \draw [<->,thick] (0,5) -- (0,0) -- (3.0,0);
			\draw [-,thick] (0,0) -- (-3.5,0);
			\draw [-,thin] (0,2.15) -- (3.5,2.15);
			\draw [-, thin] (-3.5,2.15) -- (0,2.15);
			\draw [-,thick] (0,0) to [out=0,in=-180] (2.5,2.15);
			\draw [-,thick] (-2.5,2.15) to [out=-5,in=180] (0,0);
			
			\draw (1.5,-0.2) -- (1.5, 0.2);
			\node [below] at (1.5,-0.2) {${r}^c$};
			
			\draw (-1.5,-0.2) -- (-1.5, 0.2);
			\node [below] at (-1.5,-0.2) {${-r}^c$};

			\node [right] at (3,0) {$r$};
			\node [left] at (0,2.45) {$C^+$};
			\node [above] at (-2.4,1.20) {$g(r)$};
		  \end{tikzpicture}
		  \caption{}
		  \label{ConvexConcavea}
    \end{subfigure}
    \hskip2em
    \begin{subfigure}{.45\linewidth}
        \begin{tikzpicture}[xscale=0.6,yscale=0.6]
		    \draw [<-,thick] (0,3) -- (0,-3);
			\draw [->,thick] (-5,0) -- (3.5,0);
			\draw [-,thick] (0,0) to [out=60,in=140] (1.5,1.5) to [out=-45,in=180] (3,0.0);
			
			\draw [-,thick] (-3.0,-0.0) to [out=0,in=130] (-1.5,-1.5) to [out=-50, in=245] (0,0);
			
			\draw (3.0,-0.2) -- (3.0, 0.2);
			\draw (-3.0,-0.2) -- (-3.0, 0.2);
			\draw (1.2,-0.2) -- (1.2, 0.2);
			\draw (-1.2,-0.2) -- (-1.2, 0.2);
			\node [below] at (1.2,-0.2) {${r}^c$};
			\node [below] at (-1.2,-0.2) {$-{r}^c$};
			\node [below] at (3.0,-0.2) {${r}^+$};
			\node [below] at (-3.0,-0.2) {$-{r}^+$};
			\node [right] at (3.5,0) {${r}$};
			\node [right] at (0,2.2) {$g'(r)$};
		  \end{tikzpicture}
		   \caption{}
		   \label{ConvexConcaveb}
    \end{subfigure}
    \caption{\bf (a) The potential function $g(r)$  for tensile force. Here $C^+$ is the asymptotic value of $g$. (b) Cohesive force. The derivative of the force potential goes smoothly to zero at $\pm r^+$.}\label{ConvexConcave}
\end{figure}

The force potential is a function of the strain and is defined for all $\bx,\by$ in $D$ by
\begin{equation}\label{tensilepot}
\mathcal{W}^\epsilon (S(\by,\bx,\bu(t)))=J^\epsilon(|\by-\bx|)\frac{1}{\epsilon^3\omega_2|\by-\bx|}g(\sqrt{|\by-\bx|}S(\by,\bx,\bu(t)))
\end{equation}
where $\mathcal{W}^\epsilon (S(\by,\bx,\bu(t)))$ is the pairwise force potential per unit length between two points $\bx$ and $\by$. It is described in terms of its potential function $g$, given by $g(r)=h(r^2)$ where $h$ is concave, see \autoref{ConvexConcave}. Here $\omega_2$ is the area of the unit disk and $\epsilon^2\omega_2$ is the area of the horizon $\mathcal{H}_{\epsilon}(\bx)$.

The potential function $g$ represents a convex-concave potential such that the associated force acting between material points $\bx$ and $\by$ are initially elastic and then soften and decay to zero as the strain between points increases, see \autoref{ConvexConcave}. The first well for $\mathcal{W}^\epsilon (S(\by,\bx,\bu(t)))$ is at zero tensile strain and the potential function satisfies
\begin{align}
\label{choice at oregon}
g(0)=g'(0) = 0.
\end{align}
The well for $\mathcal{W}^\epsilon (S(\by,\bx,\bu(t)))$ in the neighborhood of infinity is characterized by the horizontal asymptote $\lim_{S\rightarrow \infty} g(S)=C^+$, see \autoref{ConvexConcave}. The critical tensile strain $S_c>0$ for which the force begins to soften is given by the inflection point ${r}^c>0$ of $g$ and is 
\begin{equation}
S_c=\frac{{r}^c}{\sqrt{|\by-\bx|}},
\label{crittensileplus}
\end{equation}
and $S_+$ is the strain at which the force goes to zero 
\begin{equation}
S_+=\frac{{r}^+}{\sqrt{|\by-\bx|}}.
\label{crittensileplus2}
\end{equation}
We assume here that the  potential functions are bounded and are smooth.

\subsection{Peridynamic equation of motion}
The potential energy of the motion is given by
\begin{equation}\label{new peri}
\begin{aligned}
PD^\epsilon(\bu)=\int_D \int_{\mathcal{H}_\epsilon(\bx)\cap D} |\by-\bx|\mathcal{W}^\epsilon(S(\by,\bx,\bu(t)))\,d\by d\bx.
\end{aligned}
\end{equation}
In this treatment the material is assumed homogeneous and the density $\rho$ is constant. 
The set notation $H_\epsilon(\bx)\cap D$ means if $\bx$ belongs to $D$ and if the line connecting $\bx$ to $\by$ crosses the boundary $\partial D$ then the strain and the energy  $\mathcal{W}^\epsilon(S(\by,\bx,\bu(t)))$ associated with these two points is zero. We consider single edge notched specimen  $D$ pulled apart by a body force on the top and bottom of the domain consistent with plain strain loading. In the nonlocal setting the ``traction'' is given by an $\delta$ thick layer of body force on the top and bottom of the domain. For this case the body force is written as 
\begin{equation}\label{bodyforce}
\begin{aligned}
\bb(\bx,t)=\be^2\delta^{-1}{g_+(x_1,t)}\chi_+^\delta(x_1,x_2) \hbox{  on the top layer and}\\
\bb(\bx,t)=\be^2\delta^{-1}{g_-(x_1,t)}\chi_-^\delta(x_1,x_2)\hbox{  on the bottom  layer,}
\end{aligned}
\end{equation}
where $\be^2$ is the unit vector in the vertical direction, $\chi^\delta_+$ and $\chi^\delta_-$ are the characteristic functions of the boundary layers given by
\begin{equation}\label{layers}
\begin{aligned}
\chi_+^\delta(x_1,x_2)&=1\,\,\,\hbox{on $\{0<x_1<a, \,b/2-\delta<x_2<b/2\}$ and $0$ otherwise,}\\
\chi_-^\delta(x_1,x_2)&=1\,\,\,\hbox{on $\{0<x_1<a, \,-b/2<x_2<-b/2+\delta\}$ and $0$ otherwise},
\end{aligned}
\end{equation}
and the top and bottom traction forces are equal and in opposite directions, ie., $g_-(x_1,t)=-g_+(x_1,t)$ and
$g_+(x_1,t)>0$. We take the functions $g_-$ and $g_+$ to be smooth in the variables $x_1$  and $t$ such that
\begin{equation}\label{l2}
\sup_{t\in[0,T]}\{\Vert \bb(\bx,t)\Vert^2_{L^2(D;\mathbb{R}^2)}\}\,<\infty.
\end{equation}

For any in-plane rigid body motion $\bw(\bx)=\Omega\times\bx+\bc$ where $\Omega$ and $\bc$ are constant vectors we see that
\begin{equation}
\label{rigid}
\begin{aligned}
\int_D\bb\cdot\bw\,d\bx=0  \hbox{        and        }  S(\by,\bx,\bw)=0,
\end{aligned}
\end{equation}
and for future reference we denote the space of all square integrable fields orthogonal to rigid body motions in the $L^2$ inner product  by 
\begin{equation}
\label{rigidbody}
\dot L^2(D;\mathbb{R}^2).
\end{equation}

We define the {Lagrangian}  $${\rm{L}}(\bu,\partial_t \bu,t)=\frac{\rho}{2}||\dot \bu||^2 _{L^2 (D;\mathbb{R}^2)}-PD^\epsilon(\bu)+\int_D \bb^\epsilon\cdot \bu \,d\bx,$$
where $\dot \bu=\frac{\partial \bu}{\partial t}$ is the velocity and $\Vert \dot \bu\Vert_{L^2(D;\mathbb{R}^2)}$ denotes the $L^2$ norm of the vector field $\dot \bu: D\rightarrow \mathbb{R}^2$.
We write the action integral for a time evolution over the interval $0<t<T,$ 
\begin{equation}\label{action}
\begin{aligned}
I=\int_0^T{\rm{L}}(\bu,\partial_t \bu,t)\,dt,
\end{aligned}
\end{equation}
We suppose $\bu^\epsilon(t)$ is a stationary point and $\bw(t)$ is a perturbation and applying the {principal of least action} gives the nonlocal dynamics
\begin{eqnarray}\label{energy based weakform}
\begin{aligned}
&\rho \int_0^T\int_D\dot{\bu}^\epsilon(\bx,t)\cdot {\dot\bw}(\bx,t)d\bx\,dt\\
 &=
\int_0^T\int_D \int_{H_\epsilon(\bx)\cap D} |\by-\bx|\partial_S\mathcal{W}^\epsilon(S(\by,\bx,\bu^\epsilon(t)))S(\by,\bx,\bw(t))\,d\by d\bx\,dt\\
&-\int_0^T\int_D\bb(\bx,t)\cdot\bw(\bx,t)d\bx\,dt.
\end{aligned}
\end{eqnarray}
and an integration by parts gives the strong form 
\begin{equation}\label{energy based model2}
\begin{aligned}
\rho \ddot{\bu}^\epsilon(\bx,t)=\mathcal{L}^\epsilon(\bu^\epsilon)(\bx,t)+\bb(\bx,t),\hbox{  for  $\bx\in D$}.
\end{aligned}
\end{equation}
Here $\mathcal{L}^\epsilon(\bu^\epsilon)$ is the peridynamic force 
\begin{align}\label{pdforce}
\mathcal{L}^\epsilon(\bu^\epsilon)=\int_{\mathcal{H}_{\epsilon}(x)\cap D} {\bolds{f}}^\epsilon(\by,\bx)\;d\by
\end{align}
and ${\bolds{f}}^\epsilon(\bx,\by)$ is given by
\begin{align}\label{nonlocforcetensite}
{\bolds{f}}^\epsilon(\bx,\by)&=2\partial_S\mathcal{W}^\epsilon(S(\by,\bx,\bu^\epsilon(t)))\be_{\by-\bx},
\end{align}
where
\begin{align}\label{derivbond}
\partial_S\mathcal{W}^\epsilon(S(\by,\bx,\bu^\epsilon(t)))&=\frac{1}{\epsilon^{3} \omega_2}\frac{J^\epsilon(|\by-\bx|)}{|\by-\bx|}\partial_S g(\sqrt{|\by-\bx|}S(\by,\bx,\bu^\epsilon(t))).
\end{align}

The dynamics is complemented with the initial data 
\begin{align}\label{idata}
\bu^\epsilon(\bx,0)=\bu_0(\bx), \qquad  \partial_t \bu^\epsilon(\bx,0)=\bv_0(\bx).
\end{align}
Where $\bu_0$ and $\bv_0$ lie in $\dot L^2(D;\mathbb{R}^2)$. 

For reference we now introduce the standard space  $C^2([0,T];\dot L^2(D;\mathbb{R}^2))$  given by all the functions $\bv(\bx,t)$ for which $\bv(\bx,t)$, $\dot\bv(\bx,t)$, $\ddot\bv(\bx,t)$ belong to $\dot L^2(D;\mathbb{R}^2)$ for $0\leq t\leq T$ and
\begin{equation}\label{bounded1st}
\begin{aligned}
&\sup_{0\leq t \leq T}\Vert \bu^\epsilon(\bx,t)\Vert_{L^2(D;\mathbb{R}^2)}<\infty,\\
&\sup_{0\leq t\leq T}\Vert \dot\bu^\epsilon(\bx,t)\Vert_{L^2(D;\mathbb{R}^2)}<\infty,\\
&\sup_{0\leq t\leq T}\Vert \ddot\bu^\epsilon(\bx,t)\Vert_{L^2(D;\mathbb{R}^2)}<\infty.
\end{aligned}
\end{equation}
The initial value problem for the nonlocal evolution  given by \eqref{energy based model2} and \eqref{idata}
or equivalently by \eqref{energy based weakform} and \eqref{idata} is seen to have a unique solution in \\$C^2([0,T];\dot L^2(D;\mathbb{R}^2))$, see section \ref{existenceuniqueness}. The nonlocal evolution $\bu^\epsilon(\bx,t)$ is uniformly bounded in the mean square norm over bounded time intervals $0<t<T$, i.e., 
\begin{eqnarray}
\max_{0<t<T}\left\{\Vert \bu^\epsilon(\bx,t)\Vert_{L^2(D;\mathbb{R}^2)}^2\right\}<K,
\label{bounds}
\end{eqnarray}
where the upper bound $K$ is independent of $\epsilon$ and depends only on the initial conditions and body force applied up to time $T$. 
This follows from Theorem 2.3 of \cite{CMPer-Lipton}.

\section{Fracture toughness and elastic properties for the cohesive model: as specified through the force potential}\label{toughelastic}
\begin{figure} 
\centering
\begin{tikzpicture}[xscale=0.78,yscale=0.78]
\draw [-,dashed, thick] (0,3.5) -- (0,-0.75);
\draw [->,ultra thick] (0,.95) -- (0,0);
\draw [->,thick] (0,0) -- (2.8,.95);
\draw[->,thick] (0,0) -- (-1.7,1.7);
\draw [-,thick] (-3.9,.95) -- (3.9,.95);
\node [right] at (-1.7,1.7) {${\bf y}$};
\node [left] at (0,.55) {${z}$};

\draw [ultra thick] (2.7,1) arc [radius=3, start angle=25, end angle=155];
\node [right] at (0,-0.2) {${\bf x}$};
\draw [->,domain=90:25] plot ({2.19*cos(\x)}, {2.19*sin(\x)});

\draw [-] (0,0) -- (2.0,0.95);
\node [above] at (2.05,0.85) {$\zeta$};
\draw [->] (.9,3.5) arc (0:315:.85cm and .3cm);
\node [below] at (2.5,.75) {$\epsilon$};
\node [below] at (0.2,-0.75) {$\epsilon$};
\draw [-,thick] (-2.9,-0.75) -- (2.9,-0.75);
\node [right] at (.8,3.3) {$\theta$};
\node [right] at (-0.19,1.50) {$\scriptscriptstyle{\cos^{-1}(z/\zeta)}$};
\node [above] at (-3.5,0.95) {\footnotesize{Crack}};
\end{tikzpicture} 
\caption{{\bf Evaluation of fracture toughness $\mathcal{G}_c$. For each point $\bx$ along the dashed line, $0\leq z\leq \epsilon$, the work required to break the interaction between $\bx$ and $\by$ in the spherical cap is summed up in \eqref{calibrate2formula} using spherical coordinates centered at $\bx$.}}
 \label{compute}
\end{figure}
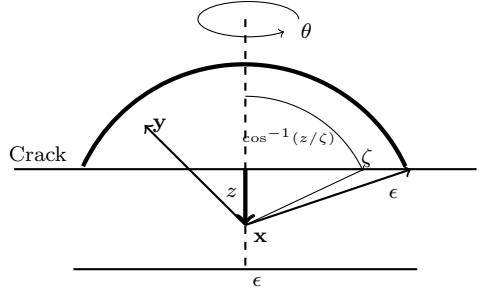

For finite horizon $\epsilon>0$ the  fracture toughness and elastic moduli are recovered 
directly from the cohesive strain potential $\mathcal{W}^\epsilon (S(\by,\bx,\bu(t)))$.
Here the fracture toughness $\mathcal{G}_c$ is defined to be  the energy per unit length required eliminate interaction between each point $\bx$ and $\by$ on either side of a line in $\mathbb{R}^2$. Because of the finite length scale of interaction only the force between pairs of points within an $\epsilon$ distance from the line are considered.  The  fracture toughness $\mathcal{G}_c$ is calculated in \cite{CMPer-Lipton}. It is given by the formula
\begin{eqnarray}
\mathcal{G}_c=2\int_0^\epsilon\int_z^\epsilon\int_0^{\arccos(z/\zeta)}\mathcal{W}^\epsilon(\mathcal{S}_+)\zeta^2\,d\psi\,d\zeta\,dz
\label{calibrate2formula}
\end{eqnarray}
where $\zeta=|\by-\bx|$, see Figure \ref{compute}. Substitution of  $\mathcal{W}^\epsilon (S(\by,\bx,\bu(t)))$ given by \eqref{tensilepot}  into \eqref{calibrate2formula} and calculation delivers the  formula 
\begin{eqnarray}
\mathcal{G}_c=\, \frac{4}{\pi}\int_{0}^1h(S_+)r^2J(r)dr.
\label{epsilonfracttough}
\end{eqnarray}
It is evident from this calculation that the fracture toughness is the same for all choice of horizons. This provides the rational behind the $\epsilon$ scaling of the potential \eqref{tensilepot} for the cohesive model. Moreover the layer width on either side of the crack centerline over which the force is applied to create new surface tends to zero with $\epsilon$. In this way $\epsilon$ can be interpreted as a parameter associated with the extent of the process zone of the material.

To find the elastic moduli associated with the cohesive potential\\ $\mathcal{W}^\epsilon (S(\by,\bx,\bu(t)))$ we
suppose the displacement inside $\mathcal{H}_\epsilon(\bx)$ is affine, that is, $u(\bx)=F\bx$ where $F$ is a constant matrix.   For small strains, i.e., $\mathcal{S}=Fe\cdot e\ll\mathcal{S}_c$, a Taylor series expansion at zero strain shows that the strain potential is linear elastic to leading order and characterized by elastic moduli $\mu$ and $\lambda$ associated with a linear elastic isotropic material  
\begin{eqnarray}
W(\bx)&=&\int_{\mathcal{H}_\epsilon({\bx})}|\by-\bx|\mathcal{W}^\epsilon(S(\by,\bx,\bu))\,d\by\nonumber\\
&=&\mu |F|^2+\frac{\lambda}{2} |Tr\{F\}|^2+O(\epsilon|F|^4).
\label{LEFMequality}
\end{eqnarray}
The elastic moduli $\lambda$ and $\mu$ are calculated directly from the strain energy density and are given by
\begin{equation}\label{calibrate1}
\mu=\lambda=M\frac{1}{4} h'(0) \,, 
\end{equation}
where the constant $M=\int_{0}^1r^{2}J(r)dr$.

\section{Convergence of the nonlocal fracture evolution to the local fracture evolution}\label{prelim}
In this section we describe the nature of the convergence of the solutions of the nonlocal cohesive model to the local fracture evolution and exhibit properties of the local fracture evolution. In the next section we will identify the interaction between the crack set and the intact material to see that the  local fracture model described in \cite{Freund} is recovered as the limit of nonlocal cohesive evolutions.  In this section we show that in the limit of small horizon $\epsilon\rightarrow 0$ the peridynamic solutions converge in the $L^2$ norm to displacements that are 
linear elastodynamic off the crack set, that is, the PDEs of the local theory hold at points off of the crack.
The elastodynamic balance laws are characterized by elastic moduli $\mu$, $\lambda$.
The evolving crack set possesses bounded Griffith surface free energy associated with the fracture toughness $\mathcal{G}_c$. The results presented below were obtained earlier in  
\cite{CMPer-Lipton3} and \cite{CMPer-Lipton}.

Suppose the initial values are taken the same for all $\epsilon$  and the initial displacement field $\bu_0(\bx)$ and initial velocity field $\bv_0(\bx)$  are in $\dot L^2(D;\mathbb{R}^2)$.  Now consider the sequence of solutions $\bu^{\epsilon}$ of the initial value problem \eqref{energy based weakform} and \eqref{idata} with progressively smaller peridynamic horizons $\epsilon$.
It is assumed as in \cite{CMPer-Lipton} that the magnitude of the displacement $\bu^\epsilon$ is bounded uniformly in $(\bx,t)$ for all horizons $\epsilon>0$.
On passing to a subsequence $\{\epsilon_n\}_{n=1}^\infty$ if necessary the  peridynamic evolutions $\bu^{\epsilon_n}$   converge in mean square uniformly in time to a 
limit evolution $\bu^0(\bx,t)$ in $C([0,T];\dot L^2(D;\mathbb{R}^2 )$ i.e.,
\begin{eqnarray}
\lim_{\epsilon_n\rightarrow 0}\max_{0\leq t\leq T}\int_D|\bu^{\epsilon_n}(\bx,t)-\bu^0(\bx,t)|^2\,d\bx=0,
\label{unifconvg}
\end{eqnarray}
moreover $\bu_t^0(\bx,t)$ belongs to $L^2(0,T;\dot L^2(D;\mathbb{R}^{2}))$.
The limit evolution $\bu^0(\bx,t)$ is also found to have bounded Griffith surface energy and elastic energy given by
\begin{eqnarray}
\int_{D}\,\mu |\mathcal{E} \bu^0(t)|^2+\frac{\lambda}{2} |{\rm div}\,\bu^0(t)|^2\,d\bx+\mathcal{G}|\mathcal{J}_{\bu^0(t)}|\leq C, 
\label{LEFMbound}
\end{eqnarray}
for $0\leq t\leq T$,
where $\mathcal{J}_{\bu^0(t)}$ denotes the evolving fracture  surface inside the domain $D$, across which the  displacement $\bu^0$ has a jump discontinuity and $|\mathcal{J}_{\bu^0(t)}|$ is its length.

In domains away from the crack set the limit evolution satisfies local linear elastodynamics (the PDEs of the standard
theory of solid mechanics).
Fix a tolerance $\tau>0$. 
If for interior subdomains $D'\subset D$ and for times $0<t<T$ the associated strains ${S}(\by,\bx,\bu^\epsilon)$ satisfy $|{S}(\by,\bx,\bu^\epsilon)|<S_c$ for every $\epsilon<\tau$  
then it is found that the limit evolution $\bu^0(\bx,t)$  is governed by the PDE
\begin{eqnarray}
\rho \bu^0_{tt}(\bx,t)= { div\,}\sigma(\bx,t), \hbox{on $[0,T]\times D'$},
\label{waveequationn}
\end{eqnarray}
where the stress tensor $\sigma$ is given by
\begin{eqnarray}
\sigma =\lambda I_2 Tr(\mathcal{E}\,\bu^0)+2\mu \mathcal{E}\bu^0,
\label{stress}
\end{eqnarray}
$I_2$ is the identity on $\mathbb{R}^2$, and $Tr(\mathcal{E}\,\bu^0)$ is the trace of the strain. The identity \eqref{waveequationn} is seen to hold in the distributional sense. The convergence of the peridynamic equation of motion to the local linear elastodynamic equation away from the crack set is 
consistent with the convergence of peridynamic equation of motion for {\emph{convex}} peridynamic potentials 
as seen in \cite{CMPer-Emmrich2}, \cite{CMPer-Silling4}, and  \cite{CMPer-Mengesha2}. 

We can say more about the limit displacement $\bu^0$ in light of the symmetry of the domain and boundary loads treated here. To refine our description we introduce the zone inside the domain where the force between two points is decreasing with increasing strain.
First we fix the time $t$ and decompose the strain $S(\by,\bx,\bu^\epsilon(t))$ as 
\begin{eqnarray}
S(\by,\bx,\bu^\epsilon(t))=S(\by,\bx,\bu^\epsilon(t))^- +S(\by,\bx,\bu^\epsilon(t))^+
\label{decomposeS}
\end{eqnarray}
where
\begin{equation}
S(\by,\bx,\bu^\epsilon(t))^-=\left\{\begin{array}{ll}S(\by,\bx,\bu^\epsilon(t)),&\hbox{if  $|S(\by,\bx,\bu^\epsilon(t))|<S_c$}\\
0,& \hbox{otherwise}
\end{array}\right.
\label{decomposedetailsS}
\end{equation}
and
\begin{equation}
S(\by,\bx,\bu^\epsilon(t))^+=\left\{\begin{array}{ll}S(\by,\bx,\bu^\epsilon(t)),&\hbox{if  $|S(\by,\bx,\bu^\epsilon(t))|\geq S_c$}\\
0,& \hbox{otherwise}
\end{array}\right.
\label{decomposedetailsS2}
\end{equation}

The subset of points $\bx\in D$ for which there is at least one $\by\in {\mathcal{H}}_{\epsilon}(\bx)\cap D$ so that $|S(\by,\bx,\bu^\epsilon(t))|\geq S_c$ is denoted by $SZ^\epsilon(t)$. This is the set of points in $D$ for which there are always points $\by$ inside $\mathcal{H}_\epsilon(\bx)$ for which the force between $\bx$ and $\by$ is decreasing with increasing strain. Motivated by the symmetry of the domain and the loading on the top and bottom boundaries we make the following geometric hypothesis on $SZ^\epsilon(t)$:
\begin{hypothesis}\label{hyp1}
Suppose $\bx$ belongs to $SZ^\epsilon$ then if $\by\in {\mathcal{H}}_{\epsilon}(\bx)\cap D$ and $\bx$ lie on different sides of the $x_2=0$ axis then $|S(\by,\bx,\bu^\epsilon(t))|\geq S_c$, on the other hand if the points $\by\in {\mathcal{H}}_{\epsilon}(\bx)\cap D$ and $\bx$ lie on the same side then $|S(\by,\bx,\bu^\epsilon(t))|<S_c$. 
 \end{hypothesis}
 This hypothesis is supported by the numerical simulations, see Figures \ref{fig:softzone} and \ref{fig:softzoneoverlap} where $SZ^\epsilon$ emerges from the simulation.
\noindent  One easily sees  that $SZ^\epsilon(t)$ is contained in a thin rectangle about $x_2=0$ of the form$\{\ell(0)\leq x_1\leq a^\epsilon,\,-\epsilon<x_2<\epsilon\}$. In light of this hypothesis one can prove the  following limit relating the sequence of strains $S(\by,\bx,\bu^\epsilon(t))$ to the jump set of the limit evolution ${\mathcal{J}}_{\bu^0(t)}$.

\begin{proposition}
\label{measures}
\begin{equation}\label{convinmeasure}
\begin{aligned}
&\lim_{\epsilon_n\rightarrow 0}\frac{1}{{\epsilon_n}^2\omega_2} \int_D\int_{{\mathcal H}_{\epsilon}(\bx)\cap D}\frac{|\by-\bx|}{\epsilon_n}J^{\epsilon_n}(|\by-\bx|)S(\by,\bx,\bu^{\epsilon_n}(t))^-d\by\,\varphi(\bx)\,d\bx
\\
&=\int_D div\,\bu^0(\bx,t)\varphi(\bx)\,d\bx\\
&\lim_{\epsilon_n\rightarrow 0}\frac{1}{{\epsilon_n}^2\omega_2} \int_{SZ^\epsilon}\int_{{\mathcal{H}}_{\epsilon}(\bx)\cap D}\frac{|\by-\bx|}{\epsilon_n}J^{\epsilon_n}(|\by-\bx|)S(\by,\bx,\bu^{\epsilon_n}(t))^+d\by\,\varphi(\bx)\,d\bx\\
 &=\int_{{\mathcal{J}}_{\bu^0(t)}}[\bu^0(\bx,t)]\cdot \bn\,\varphi(\bx)d{\mathcal{H}}^1(\bx)
\end{aligned}
\end{equation}
for all scalar test functions $\varphi$ that are differentiable with support in $D$. Here $[\bu^0(\bx,t)]$ denotes the jump in displacement.
\end{proposition}
This proposition is proved in section \ref{lengthjump}. To summarize the nonlocal evolutions converge as $\epsilon\rightarrow 0$ to a limiting evolution which has bounded Griffith surface energy and elastic energy at every time in the evolution and satisfies the wave equation away from the jump set $\mathcal{J}_{\bu^0}$. 
In the next section we exhibit the explicit interaction between the crack and the intact domain.

\section{Interaction between crack and undamaged material in the local limit}\label{ss:crackwaveinteraction}
We begin by defining a crack center line for the nonlocal model.
\begin{definition}\label{def1}
We fix the horizon size and introduce the notion of the crack center line at time $t$ for the nonlocal model. Recall that the right most boundary of the notch has the coordinates $x_1=\ell(0),x_2=0$ and the crack centerline is given by the interval $C^\epsilon(t)=\{\ell(0)\leq x_1\leq \ell^\epsilon(t), \,x_2=0\}$ across which the force between two points $\bx$ and $\by$ is zero if the line segment connecting them intersects this interval. 
\end{definition}

The crack center line is  centered in the softening zone $SZ^\epsilon$.
It is reiterated that in the nonlocal formulation the crack center line is part of the solution and its  location, shape, and evolution emerges from the numerical simulations. The numerical computation of the initial boundary value problem \eqref{energy based model2} and \eqref{idata} is given in section \ref{numerical} and provides the motivation behind hypotheses \ref{hyp1}, and the hypotheses \ref{hyp2}, and \ref{hyp3} given in this section. Motivated by the numerical experiments for this problem  we now make the following hypothesis.

\begin{figure}
  \centering
  \includegraphics[width=0.5\linewidth]{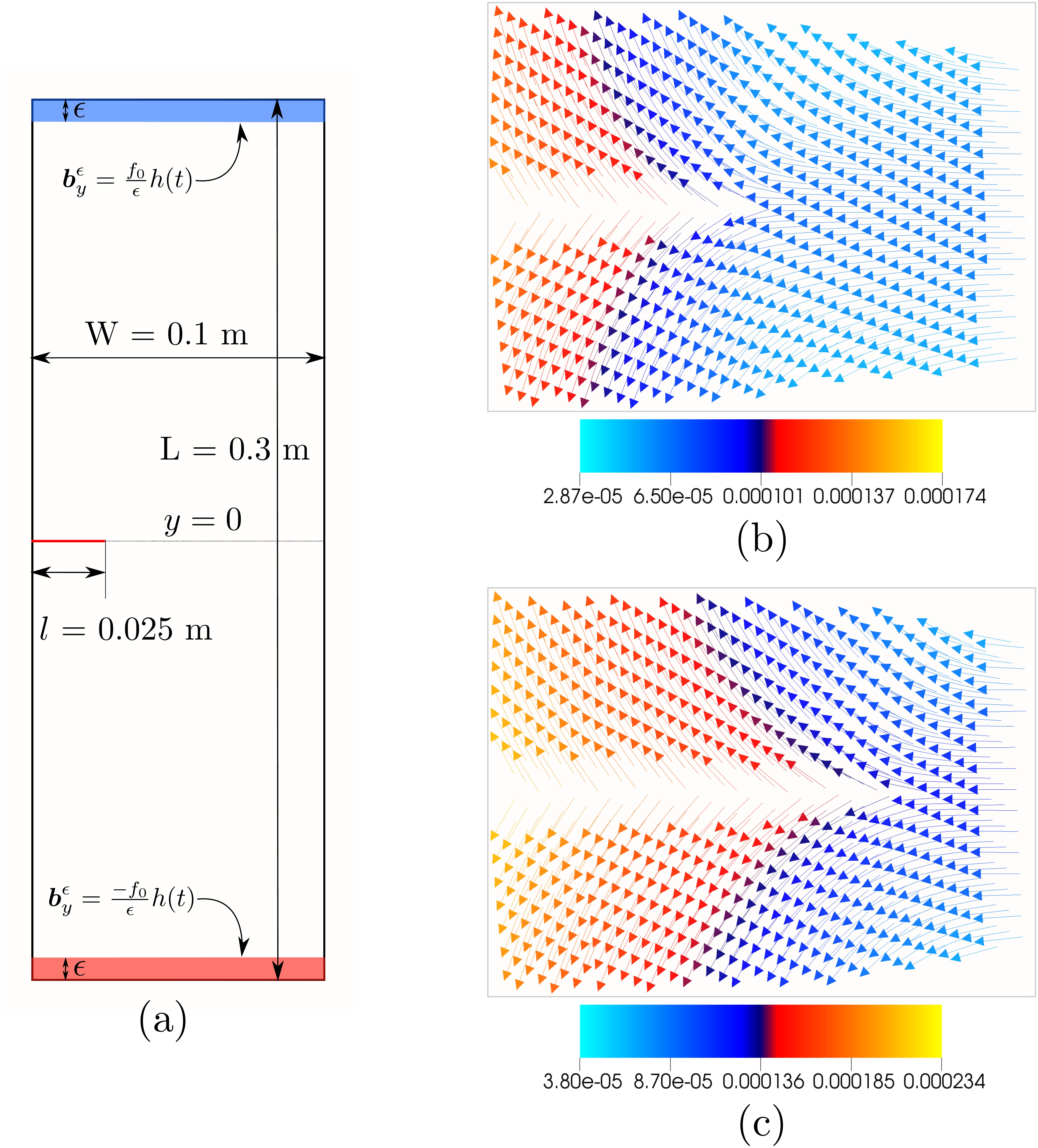}
  \caption{{\bf (a) Setup. (b), (c) Displacement profile of nodes in the domain $[0, 0.1]\times [-0.03, 0.03]$ at time $t = 460\, \mu$s and $t = 520\, \mu$s for horizon $\epsilon = 0.625$ mm}.}
  \label{fig:uplot}
\end{figure}

We define $F^\epsilon(t)$ to be the subset of points $\bx\in D$ for which there is a $\by\in {\mathcal{H}}_{\epsilon}(\bx)\cap D$ such that $S(\by,\bx,\bu^\epsilon(t))>S_+$. $F^\epsilon(t)$ is the collection of $\bx$ such that the force between at least one of its neighbors $\by$ is zero. This is called bond failure between $\by$ and $\bx$.
\begin{hypothesis}\label{hyp2}
We suppose that $SZ^\epsilon=F^\epsilon$, i.e., once bonds soften they fail.
\end{hypothesis}
We see that $SZ^\epsilon=F^\epsilon$ contains the thin rectangle $\{\ell(0)\leq x_1<\ell^\epsilon(t), \,-\epsilon<x_2<\epsilon\}$ and the width of the zone containing all pairs $\bx$ and $\by$ for which $S(\by,\bx,\bu^\epsilon(t))>S_+$ and there is negligible force between them. The width of this rectangle is $2\epsilon$. This is seen in the simulation, see Figure \ref{fig:softzone}. Thus the displacements adjacent to this zone are not influenced by forces on the other side of the zone. The physical picture is illustrated by the numerical examples  that shows that there is no difference in displacement parallel to the crack center line across the $x_2=0$ axis and a Mode I crack is propagating,  pulled open by the body loads at the top and bottom of the domain, see
figures \ref{fig:uplot}(b) and \ref{fig:uplot}(c).

This motivates the next hypothesis on the displacement and strain  adjacent to the crack center line.
\begin{hypothesis}\label{hyp3}
The displacement is always directed away from the crack center line on either side adjacent to it and across the crack centerline the strain grows as $\epsilon^{-1}$, i.e., 
\begin{equation}\label{nointerpennatration}
\frac{1}{\epsilon^2\omega_2}\int_{{\mathcal H}_\epsilon(\bx)\cap D}\frac{|\by-\bx|}{\epsilon}J^\epsilon(|\by-\bx|)S(\by,\bx,\bu^\epsilon)\,d\by>\frac{\alpha}{\epsilon}>0,
\end{equation}
for all $\bx$ in $SZ^{\epsilon/2}(t)$ and
\begin{equation}\label{nointerpennatrationlim}
[\bu^0(\bx,t)]\cdot \bn>\alpha>0,\hbox{  for all $\bx\in\mathcal{J}_{\bu^0(t)}$}.
\end{equation}
\end{hypothesis}
On the other hand the magnitude of the nonlocal strain $|S(\by,\bx,\bu^\epsilon)|$ is bounded by $r^c/\sqrt{|\by-\bx|}$ between all points $\bx$ and $\by$ not separated from each other by the crack center line.

We now begin by formally describing the convergence of the solutions $\bu^\varepsilon$ of the nonlocal  initial boundary problem \eqref{energy based weakform} and \eqref{idata} to the local evolution $\bu^0$. The limit evolution $\bu^0$
is a solution of
\begin{equation}\label{formal}
\rho \ddot\bu^0(t)=div\, \mathbb{C}\mathcal{E}\bu^0(t)+\bb\hbox{ on $D\setminus{\mathcal{J}_{\bu^0(t)}}$}
\end{equation}
and ${\mathcal J}_{\bu^0(t)}$ corresponds to the crack set given by the interval $\{0\leq x_1<\ell^0(t),\,\,x_2=0\}$. For $\bx \in {\mathcal J}_{\bu^0(t)}$ 
 we define $\mathbb{C}\mathcal{E}\bu^0(\bx)\bn^+=\lim_{\by\rightarrow\bx,\,y_2>0}\{\mathbb{C}\mathcal{E}\bu^0(\by)\bn\}$ and $\mathbb{C}\mathcal{E}\bu^0(\bx)\bn^-=\lim_{\by\rightarrow\bx,\,y_2<0}\{\mathbb{C}\mathcal{E}\bu^0(\by)\bn\}$. 
 The condition of zero traction of the top and bottom crack lips is given by
\begin{equation}\label{cracktraction}
 \mathbb{C}\mathcal{E}\bu^0\bn^+=0\hbox{  on  $\mathcal{J}_{\bu^0(t)}$} \hbox{  and  }\mathbb{C}\mathcal{E}\bu^0\bn^-=0\hbox{  on  $\mathcal{J}_{\bu^0(t)}$}
\end{equation}
where $\bn=\be^2$ on $\mathcal{J}_{\bu^0(t)}$. The natural boundary conditions are
\begin{equation}\label{tractions}
\begin{aligned}
 \mathbb{C}\mathcal{E}\bu^0\bn&=0\hbox{  on   $\partial D$}
 \end{aligned}
\end{equation}
where $\bn$ is the outward directed unit normal on the boundary $\partial D$.

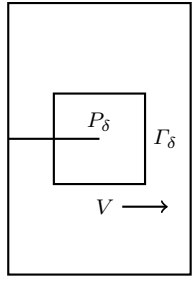
\begin{figure} 
\centering
\begin{tikzpicture}[xscale=0.60,yscale=0.60]

\draw [-,thick] (-2,0) -- (-0.0,0);

\draw [thick] (-2,-3) rectangle (2,3);

\draw [thick] (-1,-1) rectangle (1,1);

\node [above] at (0.0,0.0) {$P_\delta$};

\node [right] at (1,0.0) {$\Gamma_\delta$};

\draw [->,thick] (.5,-1.5) -- (1.5,-1.5);

\node [left] at (0.5,-1.5) {$V$};

\end{tikzpicture} 
\caption{{ \bf Contour $\Gamma_\delta$ surrounding the domain $P_\delta$ moving with the velocity $V$ of the crack centerline.}}
 \label{P}
\end{figure}

We now address  crack tip interaction with intact elastic material by calculating the change in internal energy inside a neighborhood enclosing the crack tip. Fix  a contour $\Gamma_\delta$ of diameter $\delta$  surrounding the domain $\mathcal{P}_\delta(t)$ containing the crack tip for the local model, see Figure \ref{P}. This domain is moving to the right with the crack tip velocity $V$. Proposition \ref{limitspowerflow} shows that after passing to the $\epsilon\rightarrow 0$ limit in the nonlocal cohesive model  the rate of change of internal energy inside $P_\delta$  is given by
\begin{equation}\label{energyflowbalance3bis}
\begin{aligned}
\lim_{\epsilon\rightarrow 0}\frac{d}{dt}\int_{\mathcal{P}_\delta(t)}\,T^\epsilon +W^\epsilon\,d\bx=\int_{\Gamma_\delta}\mathbb{C}\mathcal{E}\bu^0\bn\cdot\dot\bu^0\,ds-\mathcal{G}_cV(t)+O(\delta),
 \end{aligned}
\end{equation}
where $T^\epsilon$ is the kinetic energy density
and $W^\epsilon$ is the energy density given by 
\begin{equation}\label{densitynon}
W^{\epsilon}(\bx,t)=\int_{{\mathcal{H}_\epsilon(\bx)}\cap D}|\by-\bx|\mathcal{W}^{\epsilon}(S(\by,\bx,\bu^\epsilon(t)))\,d\by.
\end{equation}

If the limit of the internal energy is not changing in time then in the $\delta=0$ limit we recover the power balance for the local model given by
\begin{equation}\label{ratebalance}
\begin{aligned}
V(t) \mathcal{G}_c =J=\lim_{\delta\rightarrow 0}{J}_\delta,
 \end{aligned}
\end{equation}
where  $J_\delta$ 
is the rate of energy flowing into $\mathcal{P}_\delta(t)$ towards the crack tip for the local model given by
\begin{equation}\label{energyrate}
\begin{aligned}
 J_\delta=\int_{\Gamma_\delta}\mathbb{C}\mathcal{E}\bu^0\bn\cdot\dot\bu^0\,ds.
 \end{aligned}
\end{equation}
This well known situation can happen for cracks propagating at constant velocity inside domains with remote boundaries \cite{FreundClifton}.
Here $\bn$ is the outward directed unit normal, $ds$ is an element of arc length and $\mathbb{C}\mathcal{E}\bu^0:\mathcal{E}\bu^0$ is the elastic energy density. Clearly the crack advances if $J>0$, see Proposition \ref{limitspowerflow} below.   The rate of energy  flowing into the crack for the cohesive model is derived in section \ref{s:Jintegral}. The formula recovered here is in the form given by \cite{FreundClifton}. Based on the nature of the dynamic stress field \cite{AtaksonEshelby}, \cite{KostrovNitkin}, \\ \cite{Freund1972a}, and \cite{Willis1975} provide a representation of the instantaneous rate of energy flowing into the crack tip, is given by
\begin{equation}\label{energyratelim}
\begin{aligned}
J=\frac{1+\nu}{E}\frac{V^3}{c_s^2D}\alpha_tK_I^2(t), \end{aligned}
\end{equation}
were $\nu$ is the Poisson ratio, $E$ is the Young modulus $V$ is the crack velocity, $c_s$ is the shear wave speed, $c_l=(\lambda+2\mu/\rho)^{1/2}$ is the longitudinal wave speed, $D=4\alpha_s\alpha_l-(1+\alpha_s^2)^2$, and $\alpha_s=(1-V^2/c_s^2)^{1/2}$, $\alpha_l=(1-V^2/c_l^2)^{1/2}$.
Here $K_I(t)$ is the mode $I$ dynamic stress intensity factor and depends on the details of the loading and is not explicit. 
We arrive at the kinetic relation for the crack tip velocity given by
\begin{equation}\label{kinetic}
\begin{aligned}
\mathcal{G}_c=\frac{J}{V}=\frac{1+\nu}{E}\frac{V^2}{c_s^2D}\alpha_tK_I^2(t), \end{aligned}
\end{equation}
see \cite{Freund}. To summarize, the limit  field $\bu^0(\bx,t)$ and crack tip velocity $V(t)\be^1$ are coupled and satisfy the equations \eqref{formal}, \eqref{cracktraction}, \eqref{tractions}, and \eqref{kinetic}. Conversely these equations provide conditions necessary to determine $\bu^0(\bx,t)$ and $V(t)$.

We now follow up with the rigorous statement of the convergence of the nonlocal dynamics to the weak form of the limit problem described by \eqref{formal}, \eqref{cracktraction}, and \eqref{tractions} in Propositions \ref{convergences}, \ref{momentumlim}, and \ref{traction2}. Then in Proposition \ref{limitspowerflow} we provide the limiting balance of energy flow rate into the crack tip \eqref{ratebalance} using only the nonlocal cohesive dynamics associated with the double well potential \eqref{energy based model2}.  

Recall the crack centerline associated with a nonlocal evolution with horizon $\epsilon$ is denoted by $C^\epsilon(t)$ see definition \ref{def1}. Passing to a subsequence if necessary  we send $\epsilon$ to zero and arrive at  a limit $C^0(t)$.
We start with the rigorous proposition showing there is one unique limit of the crack center line set $C^0(t)$ and it is related to the jump set ${\mathcal J}_{\bu^0(t)}$. In what follows $|\cdot|$ denotes the one dimensional Lebesgue measure of a set and $\Delta$ is the symmetric dfference between two sets.
\begin{proposition}
\label{cracksetjumpset}
Suppose hypothesis 1 through 3 hold, then the limit of the sets $C^{\epsilon_n}=\{\ell(0)\leq x_1\leq \ell^{\epsilon_n}(t)\,,\,x_2=0\}$ exists and is given by $C^0(t)=\{\ell(0)\leq x_1< \ell^0(t)\,,\,x_2=0\}$ where 
\begin{equation}\label{cracklength}
\begin{aligned}
\ell^0(t)&=\lim_{\epsilon_n\rightarrow 0}\ell^\epsilon(t) =|{\mathcal J}_{\bu^0(t)}|\hbox{  and}\\
&|{\mathcal J}_{\bu^0(t)}\,\Delta\, C^0(t)|=0.
\end{aligned}
\end{equation}
\end{proposition}

Next we describe the convergence in terms of the suitable Hilbert space topology. The space of strongly measurable functions $\bw\,:\,[0,T]\rightarrow  \dot L^2(D;\mathbb{R}^2)$ that are square integrable in time is denoted by $L^2(0,T;\dot L^2(D;\mathbb{R}^2))$. Additionally we recall the Sobolev space $H^1(D;\mathbb{R}^2)$ with norm 
\begin{equation}\label{h1}
\Vert\bw\Vert_{H^1(D;\mathbb{R}^2)}:=\left(\int_ D\,|\bw|^2\,d\bx+\int_D|\nabla\bw|^2\,d\bx\right)^2.
\end{equation}
The subspace of $H^1(D;\mathbb{R}^2)$
containing all vector fields orthogonal to the rigid motions with respect to the $L^2(D;\mathbb{R}^2)$ inner product is written
\begin{equation}\label{h85}
\dot H^1(D;\mathbb{R}^2).
\end{equation}
The Hilbert space dual to $\dot H^1(D;\mathbb{R}^2)$ is denoted by $\dot H^1(D;\mathbb{R}^2)'$.  
The set of strongly integrable functions taking values in $\dot H^1(D;\mathbb{R}^2)'$ for $0\leq t\leq T$ is denoted by $L^2(0,T;\dot H^1(D;\mathbb{R}^2)')$.  These Hilbert spaces are well known and related to the wave equation, see \cite{Evans}. 
We pass to subsequences as necessary and the convergence is given by 
\begin{proposition}\label{convergences}
\begin{equation}\label{convgences}
\begin{aligned}
\bu^{\epsilon_n}\rightarrow \bu^0 & \hbox{  \rm strong in } C([0,T];\dot L^2(D;\mathbb{R}^2))\\
\dot\bu^{\epsilon_n}\rightharpoonup \dot\bu^0 &\hbox{ \rm weakly in  }L^2(0,T;\dot L^2(D;\mathbb{R}^2))\\
\ddot\bu^{\epsilon_n}\rightharpoonup \ddot\bu^0 &\hbox{ \rm weakly in }L^2(0,T;\dot H^1(D;\mathbb{R}^2)').
\end{aligned}
\end{equation}
\end{proposition}

The weak form of the momentum balance and zero traction condition are given are terms of a pair of variational identities over properly chosen test spaces.  For completeness we now recall the support of a vector function $\bw$  defined as $$supp\{\bw\}=\hbox{ Closure of $\{\bx\in D;\,|\bw(\bx)|\not=0\}$},$$
and we define the test space $H^{1,0}(D\setminus C^0(t),\mathbb{R}^2)\subset H^1(D;\mathbb{R}^2)$ to be the $H^1$ norm closure of the set of infinitely differentiable functions on $D$ with support sets that do not intersect the crack.
\begin{proposition}\label{momentumlim}
Suppose hypotheses 1, 2, and 3 hold then:
the field $\ddot\bu^0(t)$ is in fact a bounded linear functional on the space  $H^{1,0}(D\setminus C^0(t),\mathbb{R}^2)$ for  a.e. $t\in [0,T]$ and we have 
\begin{equation}\label{momentumlimit}
\begin{aligned}
&\rho\langle\ddot\bu^0,\bw\rangle\\
&=-\int_{D\setminus C^0(t) }\mathbb{C}\mathcal{E}\bu^0:\mathcal{E}\bw\,dx+\langle \bb,\bw\rangle,\hbox{ for all $\bw \in H^{1,0}(D\setminus C^0(t),\mathbb{R}^2)$},
\end{aligned}
\end{equation}
where $\langle\cdot,\cdot\rangle$ is the duality paring between $H^{1,0}(D\setminus C^0(t),\mathbb{R}^2)$ and its Hilbert space dual $H^{1,0}(D\setminus C^0(t),\mathbb{R}^2)'$. 
Equivalently
\begin{equation}\label{boundary tracesA}
\begin{aligned}
\mathbb{C}\,\mathcal{E}\bu^0\bn=0 &\hbox{ on $\partial D$},\\
\end{aligned}
\end{equation}
where equalities hold as elements of $H^{-1/2}(\partial D)$, and  $$\rho\ddot\bu^0=div\left(\mathbb{C}\mathcal{E}\bu^0\right)+\bb$$ as elements of $H^{1,0}(D\setminus C^0(t),\mathbb{R}^2)'$ . 
\end{proposition}
\noindent Equation \eqref{momentumlimit} is the weak formulation of the balance of momentum \eqref{formal} and \eqref{tractions}. 

We now express the weak formulation of the traction free condition at the upper and lower crack lips \eqref{cracktraction}. 
For any $\beta>0$ chosen such that $\ell(0)<\ell^0(t)-\beta$, we introduce the rectangular sets $Q^+_\beta(t)=\{\ell(0)\leq x_1<\ell^0(t)-\beta;\, 0<x_2<b/2-\delta\}$ and 
$Q^-_\beta(t)=\{\ell(0)\leq x_1<\ell^0(t)-\beta;\, -b/2+\delta<x_2<0\}$.
The Sobolev space given by the $H^1$ norm closure of the set of infinitely differentiable functions on $Q^+_\beta(t)$ with support sets that do not intersect the sets $\{\ell(0)\leq x_1<\ell^0(t)-\beta;\,x_2=b/2-\delta\}$, $\{x_1=\ell^0(t)-\beta;\, 0<x_2<b/2-\delta\}$, and $\{x_1=\ell(0);\, 0<x_2<b/2-\delta\}$ is denoted by $H^{1,0}(Q^+_\beta(t),\mathbb{R}^2)$. Similarly the $H^1$ norm closure of the set of continuously differentiable functions on $Q^-_\beta(t)$ with support sets that do not intersect the sets $\{\ell(0)\leq x_1<\ell^0(t)-\beta;\,x_2=-b/2+\delta\}$,  $\{x_1=\ell^0(t)-\beta;\, -b/2+\delta<x_2<0\}$, and  $\{x_1=\ell(0);\, -b/2+\delta<x_2<0\}$ is denoted by $H^{1,0}(Q^-_\beta(t),\mathbb{R}^2)$. The relevant trace spaces for $H^1(Q^{\pm}_\beta(t),\mathbb{R}^2)$ on $\{0\leq x_1\leq\ell^0(t)-\beta;\, x_2=0\}$ is the classic boundary space $H^{1/2}$ and its dual $H^{-1/2}$; these spaces are denoted by $H^{1/2}_\beta$ and  $H^{-1/2}_\beta$ respectfully.
\begin{proposition}\label{traction2}
Suppose hypotheses 1, 2, and 3 hold then:
the field $\ddot\bu^0(t)$ is in fact a bounded linear functional on the spaces  \,$H^{1,0}(Q^{\pm}_\beta(t),\mathbb{R}^2)$ for a.e. $t\in [0,T]$ and we have 
\begin{equation}\label{momentumlimit2}
\begin{aligned}
\rho\langle\ddot\bu^0,\bw\rangle=-\int_{D\setminus C^0(t) }\mathbb{C}\mathcal{E}\bu^0:\mathcal{E}\bw\,\,d\bx, \hbox{ for all $\bw \in H^{1,0}(Q^{\pm}_\beta(t),\mathbb{R}^2)$},
\end{aligned}
\end{equation}
where $\langle\cdot,\cdot\rangle$ is the duality paring between \,$H^{1,0}(Q^{\pm}_\beta(t),\mathbb{R}^2)$ and their duals 
\\$H^{1,0}(Q^{\pm}_\beta(t),\mathbb{R}^2)'.$  Equivalently
\begin{equation}\label{boundary traces}
\begin{aligned}
\mathbb{C}\,\mathcal{E}\bu^0\bn^+=0 \hbox{ and   }\, \mathbb{C}\,\mathcal{E}\bu^0\bn^-=0 \hbox{ on }\{\ell(0)\leq x_1<\ell^0(t)-\beta;\, x_2=0\}
\end{aligned}
\end{equation}
for all $0<\beta<\ell^0(t)-\ell(0)$ as elements of  $H^{-1/2}_\beta$ and  $$\rho\ddot\bu^0=div\left(\mathbb{C}\mathcal{E}\bu^0\right)$$ as elements of $H^{1,0}(Q^{\pm}_\beta(t),\mathbb{R}^2)'$.
\end{proposition}
\noindent Equation \eqref{boundary traces} is the weak formulation of zero traction on the crack lips \eqref{cracktraction}. 




We now investigate the power balance in regions containing the crack tip. For a nonlocal evolution with perydynamic horizon $\epsilon$ consider the rectangular contour $\Gamma^\epsilon_\delta(t)$ of diameter $\delta$ surrounding the domain $\mathcal{P}^\epsilon_\delta(t)$ containing the crack tip, see Figure \ref{Pepsilon}. We suppose $\mathcal{P}^\epsilon_\delta(t)$ is moving with the crack tip velocity $V^\epsilon(t)\be^1=\dot\ell^\epsilon(t)\be^1$ where the unit vector $\be^1$ is along the horizontal axis. Define the kinetic energy density by $T^\epsilon=\rho|\bu^\epsilon(\bx,t)|^2/2$ and the nonlocal potential energy density is given by $W^\epsilon(\bx)=\int_{\mathcal{H}_\epsilon(\bx)}|\by-\bx|\mathcal{W}^\epsilon (S(\by,\bx,\bu^\epsilon(t)))\,d\by$.

The the rate of change of internal energy inside the domain containing the crack tip  is given by
\begin{proposition}\label{changinternal}
\begin{equation}\label{powerbfornonloc}
\begin{aligned}
\frac{d}{dt}\int_{\mathcal{P}^\epsilon_\delta(t)}\,T^\epsilon +W^\epsilon\,d\bx=I^\epsilon(\Gamma^\varepsilon_\delta(t))
\end{aligned}
\end{equation}
with
\begin{equation}\label{energyflownonloc}
\begin{aligned}
I^\epsilon(\Gamma^\varepsilon_\delta(t))
=\int_{\Gamma^\epsilon_\delta(t)}\,(T^{\epsilon_n}+W^{\epsilon_n})V^{\epsilon_n}\be^1\cdot\bn\,ds-E^{\epsilon_n}(\Gamma^{\epsilon_n}_\delta(t)),
\end{aligned}
\end{equation}
and
\begin{equation}\nonumber
\begin{aligned}
&E^{\epsilon_n}(\Gamma^{\epsilon_n}_\delta(t))\\
&=\int_{A^{\epsilon_n}_\delta(t)}\int_{\mathcal{H}_{\epsilon_n}(\bx)\cap \mathcal{P}^{\epsilon_n}_\delta(t)}\partial_S\mathcal{W}^{\epsilon_n}(S(\by,\bx,\bu^{\epsilon_n}))\be_{\by-\bx}\cdot(\dot\bu^{\epsilon_n}(\bx)+\dot\bu^{\epsilon_n}(\by))\,d\by d\bx,
\end{aligned}
\end{equation}
where $\bn$ is the unit normal pointing out of the domain $\mathcal{P}^{\epsilon_n}_\delta(t)$ and $A^{\epsilon_n}_\delta(t)$  is the part of $D$ exterior to $\mathcal{P}^{\epsilon_n}_\delta(t)$.
\end{proposition}
\noindent This Proposition is derived in section \ref{s:energyrate}.

We now state the convergence of the rate of change in internal energy in the $\epsilon_n=0$ limit. 
\begin{proposition}\label{limitspowerflow}
In addition to hypothesis 1 through 3 we suppose that  $\bu^{\epsilon_n}(t)$, $\dot\bu^{\epsilon_n}(t)$, $S(\by,\bx,\bu^{\epsilon_n})$  converge uniformly to $\bu^0(t)$,  $\dot\bu^0$(t), and $\mathcal{E}\bu^0\be\cdot\be$  on subsets  away from the the crack tip for $t\in[0,T]$.  Then for  $V^{\epsilon_n}\be^1\rightarrow V\be^1$,
\begin{equation}\label{nonloctolocrate1}
\begin{aligned}
 \lim_{{\epsilon_n}\rightarrow 0}\int_{\Gamma^\epsilon_\delta(t)}\,(T^{\epsilon_n}+W^{\epsilon_n})V^{\epsilon_n}\be^1\cdot\bn\,ds&=-\mathcal{G}_c V(t)+O(\delta),\\
 \lim_{{\epsilon_n}\rightarrow 0}E^{\epsilon_n}(\Gamma^{\epsilon_n}_\delta(t))&=-\int_{\Gamma_\delta}\mathbb{C}\mathcal{E}\bu^0\bn\cdot\dot\bu^0\,ds.
 \end{aligned}
\end{equation}
The rate of change of internal energy inside the domain containing the crack tip  is given by
\begin{equation}\label{energyflowbalance3}
\begin{aligned}
\lim_{\epsilon_n\rightarrow 0}\frac{d}{dt}\int_{\mathcal{P}^{\epsilon_n}_\delta(t)}\,T^{\epsilon_n} +W^{\epsilon_n}\,d\bx=\int_{\Gamma_\delta}\mathbb{C}\mathcal{E}\bu^0\bn\cdot\dot\bu^0\,ds-\mathcal{G}_cV(t)+O(\delta).
 \end{aligned}
\end{equation}
\end{proposition}
From this proposition we see that when the rate of change of internal energy is zero the kinetic relation for  the crack tip velocity is
\begin{equation}\label{energyflowbalancedel}
\begin{aligned}
\mathcal{G}_cV(t) =\int_{\Gamma_\delta}\mathbb{C}\mathcal{E}\bu^0\bn\cdot\dot\bu^0\,ds+O(\delta),
 \end{aligned}
\end{equation}
here $\mathbb{C}\mathcal{E}\bu^0\dot\bu^0$ is energy flux into $P_\delta$.  For $J=\lim_{\delta\rightarrow 0}\int_{\Gamma_\delta}\mathbb{C}\mathcal{E}\bu^0\bn\cdot\dot\bu^0\,ds$,
we get
\begin{equation}\label{energyflowbalance}
\begin{aligned}
\mathcal{G}_cV(t)=J.
 \end{aligned}
\end{equation}
Clearly the crack moves if the rate of change of total energy of the region surrounding the crack is zero and the energy flowing into the crack tip is positive. The  semi explicit kinetic relation relating the energy flux into the crack tip and the crack velocity follows from \eqref{energyflowbalance} and is the well known one  \cite{Freund}  given by \eqref{kinetic}.

In distinction from the classic approach it is seen that  \eqref{energyflowbalancedel}, \eqref{energyflowbalance} are not postulated but instead recovered directly from \eqref{energy based model2} and are a consequence of the nonlocal cohesive dynamics in the $\epsilon_n=0$ limit.  The recovery is possible since the nonlocal model is well defined over ``the process zone'' around the crack center line tip contained inside $P^{\epsilon_n}_\delta$. The rate of change in energy internal  to $P_\delta^{\epsilon_n}$ given by \eqref{energyflownonloc}, 
\eqref{energyflowbalance3}, and the kinetic relation \eqref{energyflowbalance}
are established in section \ref{s:Jintegral}. The semi explicit form \eqref{kinetic} follows on substituting the formula for $J$ given by \cite{FreundClifton} into \eqref{energyflowbalance}.

\section{Numerical motivation for the hypotheses}
\label{numerical}
In this section we present numerical results which support the hypothesis assumed in theoretical analysis. We consider a material with density $\rho = 1200 \,kg/m^3$, Young's modulus $E = 3.24 \, GPa$, and critical energy release rate $G_c = 500 \,Jm^{-2}$. Since, this is a bond-based Peridynamic model, we are restricted to Poisson ratio $\nu = 0.25$. The pairwise interaction potential $g$ is of the form $g(r) = c (1-\exp [- \beta r^2])$. The influence function is $J(r) = 1-r$. Using equations \eqref{epsilonfracttough} and \eqref{calibrate1}, we show that $ c = 392.7$, $\beta = 1.3201 \times 10^7$. The inflection point is $r^c = 1/\sqrt{2\beta}$ and the critical strain is $S_c(\by, \bx) = r_c/\sqrt{|\by - \bx|}$. In the numerical simulation  $0<g' <<1$ at $10\times r^c$ and  $S^+(\by, \bx) = 10 \times S_c(\by, \bx)$. 

We consider the central difference time discretization on a uniform square mesh. Let $\bu^\epsilon_h$ is the approximate displacement defined on the mesh, and is given by the piecewise constant extension of the nodal displacements. The nonlocal force at mesh node $\bx_i$ is written as 
\begin{equation}
	\mathcal{L}^\epsilon(\bu^\epsilon_h)(\bx_i) = \int_{\mathcal{H}_\epsilon(\bx_i) \cap D} w(\by, \bx_i) d\by,
\end{equation}
and we approximate the force by
\begin{equation}
	\mathcal{L}^\epsilon_h(\bu^\epsilon_h)(\bx_i) \approx \sum_{\bx_j \in \mathcal{H}_\epsilon(\bx_i) \cap D} w(\bx_j, \bx_i) V_j \bar{V}_{ij},
\end{equation}
where $V_j = h^2$ is the volume represented by node $\bx_j$. The area correction $\bar{V}_{ij}$ is the ratio of the part of the area $V_j$ inside the horizon of $\bx_i$ and the area $V_j$.  In what follows the crack center line is part of the solution of the nonlocal model and the crack center line location, shape, and evolution emerges from our numerical computation of the initial boundary value problem \eqref{energy based model2} and \eqref{idata}. It is these simulations that give us the  numerical confirmation to suggest hypotheses \ref{hyp1}, \ref{hyp2}, and \ref{hyp3}.

\subsection{Crack propagation and softening zone.}
We consider a 2-d material domain $R = [0, 0.1] \times [-0.15, 0.15]$ m$^2$. We introduce a horizontal  rectangular notch of length $l = 0.025$ m originating from $(0,0)$ of width $2\epsilon$.  The time interval for simulation is $[0, 560\, \mu \text{s}]$ and the size of the time step is $\Delta t = 0.02\, \mu$s. The body force of the form $\bb(\bx, t) = (0, f_0 h(t)/\epsilon)$ and $\bb(\bx, t) = (0, -f_0 h(t))/\epsilon)$ is applied on the top and bottom layer of thickness $\delta=\epsilon$ respectively, see \autoref{fig:uplot}(a). Here $h(t)$ is a step function such that $h(t) = t$ for $t\in [0, 350 \, \mu\text{s}]$ and $h(t) = 1$ for $t\in [350 \, \mu\text{s}, 560 \, \mu\text{s}]$. We set $f_0 = 1.0\times 10^{10}$.
In the subsequent simulations we consider three different horizons, $\epsilon = 2.5, 1.25, 0.625$ mm. The mesh size for each horizon is fixed by the relation $h = \epsilon / 4$. 

\begin{figure}
  \centering
  \includegraphics[width=0.4\linewidth]{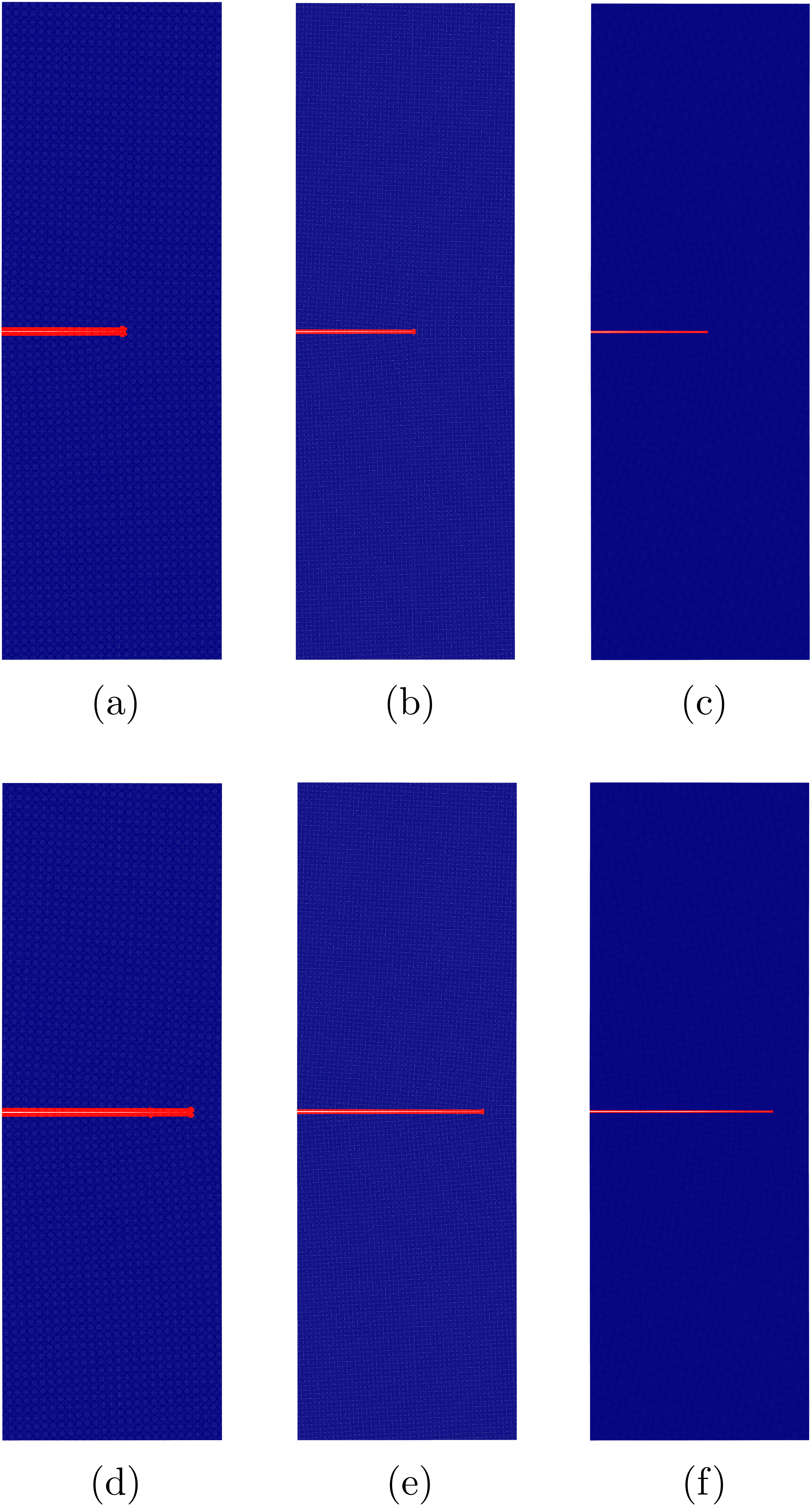}
  \caption{{\bf Softening zone (red) for different horizons. (a), (b), (c) correspond to $SZ^\epsilon(t)$ at $t=460\, \mu$s for $\epsilon = 2.5, 1.25, 0.625$ mm. (d), (e), (f) correspond to $SZ^\epsilon(t)$ at $t=520\, \mu$s for $\epsilon = 2.5, 1.25, 0.625$ mm}.}
  \label{fig:softzone}
\end{figure}

\begin{figure}
  \centering
  \includegraphics[width=0.4\linewidth]{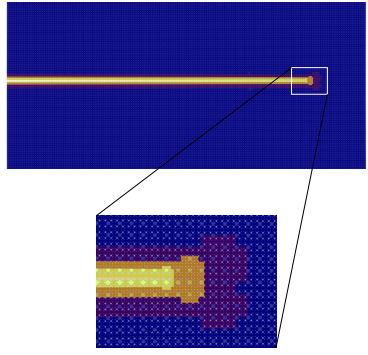}
  \caption{{\bf Top: Softening zone $SZ^\epsilon(t)$ for $\epsilon = 2.5, 1.25, 0.625\,$mm at time $t=520\,\mu$s on top of each other. Red, light yellow, and light blue color is used for $SZ^\epsilon$ of horizon $2.5, 1.25, 0.625\,$mm respectively. Bottom: Zoomed in near the crack center line tip}.}
  \label{fig:softzoneoverlap}
\end{figure}

\begin{figure}
    \centering
    \includegraphics[width=0.50\linewidth]{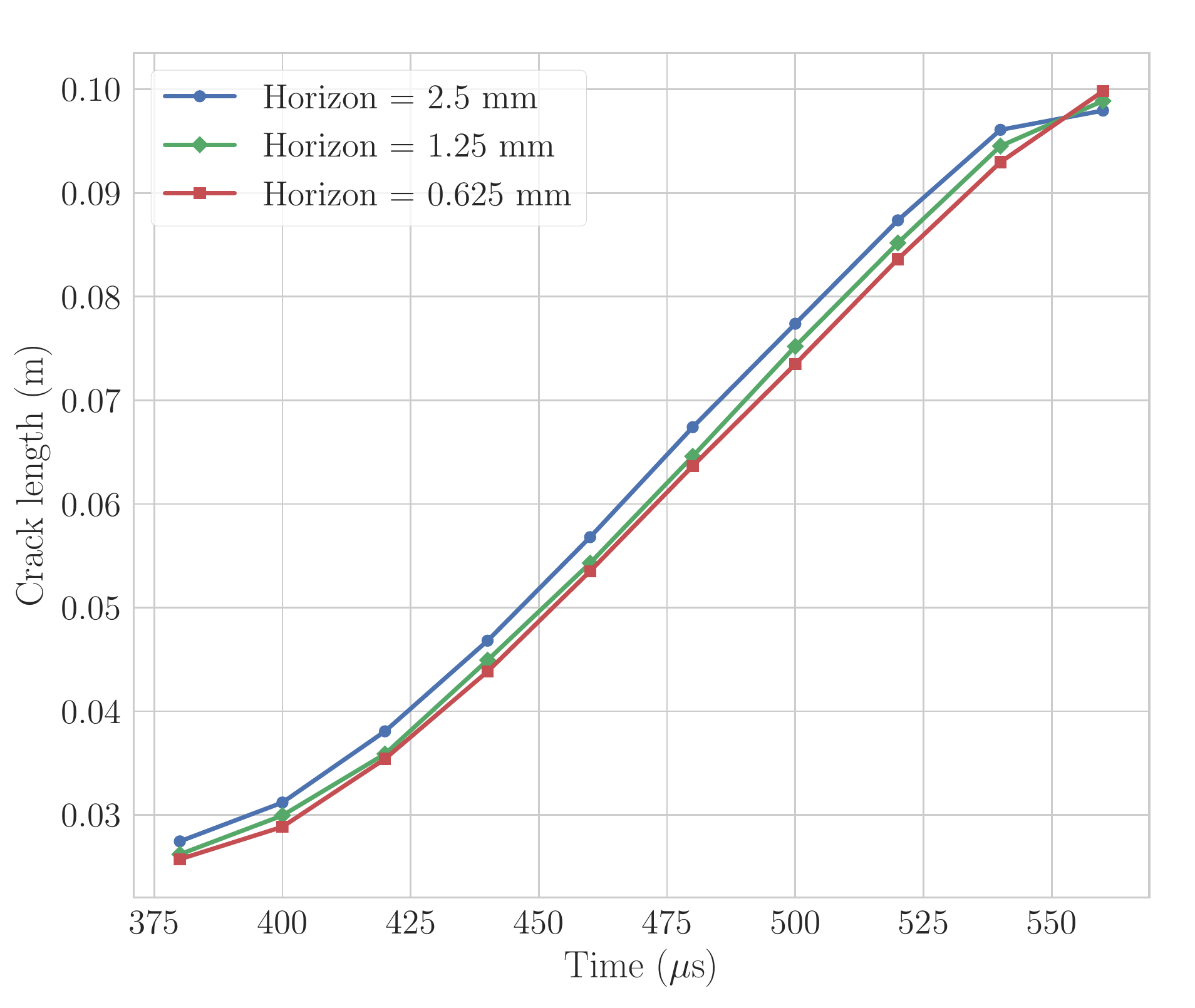}
  \caption{{\bf The crack center line length is plotted as a function of time for three different horizons.}} 
  \label{fig:cracktipalla}
\end{figure}

We use the simulations to examine the displacement profile around the crack center line.
In \autoref{fig:uplot}(b),(c) we display the displacement vector around the crack center line at times $t=460, 520\, \mu$s. From these plots, we see that nodal displacements are always pointing away from the crack line and that there is negligible force acting across the crack center line.  This provides the visual verification of Hypothesis \ref{hyp2} and \ref{hyp3}.
We now examine the softening zones that emerge from the simulations.
In \autoref{fig:softzone}, we plot the softening zone $SZ^\epsilon$ (colored in red) for three horizons at two times $460\, \mu$s and $520\, \mu$s. As one would expect, the softening zone localizes and its thickness decreases as horizon gets smaller and these results support Hypothesis \ref{hyp1}. 
We also find that the $SZ^\epsilon$ are monotone decreasing with $\epsilon$, i.e. $SZ^{\epsilon_1}(t) \subset SZ^{\epsilon_2}(t)$ when $\epsilon_1 < \epsilon_2$, see \autoref{fig:softzoneoverlap}.
The crack center line tip location  emerges from our simulations of \eqref{energy based model2}  and
in Figure \ref{fig:cracktipalla} we plot the crack center line tip location (x-coordinate) at different times. Since, $SZ^{\epsilon_1} \subset SZ^{\epsilon_2}$ when $\epsilon_1 < \epsilon_2$, we see that crack center line tip for the larger horizon is consistently ahead of the crack tip for the smaller horizon. 

\section{Existence and uniqueness of the nonlocal evolution}\label{existenceuniqueness}
We assert the existence and uniqueness for a solution $\bu^\epsilon(\bx,t)$ of the nonlocal evolution with the balance of momentum given in strong form \eqref{energy based model2}.

\begin{proposition}\label{existenceuniquness}{
\bf Existence and uniqueness of the nonlocal evolution.}
The initial value problem given by \eqref{energy based model2} and \eqref{idata} has a unique solution $\bu(\bx,t)$ such that for every $t\in [0,T]$, $\bu$ takes values in $\dot L^2(D;\mathbb{R}^2)$ and belongs to the space $C^2([0,T];\dot L^2(D;\mathbb{R}^2))$. 
\end{proposition}

The proof of this proposition follows from the Lipschitz continuity of\\ $\mathcal{L}^\epsilon(\bu^\epsilon)(\bx,t)+\bb(\bx,t)$ as a function of  $\bu^\epsilon$ with respect to the $L^2(D;\mathbb{R}^2)$ norm and the Banach fixed point theorem, see e.g. \cite{CMPer-Lipton4}.

\noindent It is remarked that in the context of the cohesive model the crack center line and $SZ^\epsilon$ describe an unstable  phase of the material. However because the peridynamic force is a Lipschitz function on $\dot L^2(D;\mathbb{R}^2)$  the model can be viewed as an ODE for vectors in $\dot L^2(D;\mathbb{R}^2)$ and is well posed.

\section{Relation between the crack set and the jump set}\label{lengthjump}

In this section we establish the relation between the fracture set and the jump set given by Proposition \ref{cracksetjumpset}. To do this we first prove Proposition \ref{measures}. We will then use this proposition together with Hypotheses 1, 2, and 3 to get Proposition \ref{cracksetjumpset}.  We start by describing the Banach space that the $\epsilon=0$ limit displacement $\bu^0$  belongs to, see \cite{CMPer-Lipton}.
The limiting $\epsilon_n\rightarrow 0$ dynamics  is given by a displacement that belongs to the space of functions of bounded deformation $SBD$.  Functions $\bu\in SBD$ belong to  $L^1(D;\mathbb{R}^d)$ (where $d=2$ in this work)  and are approximately continuous, i.e., have Lebesgue limits for almost every $\bx\in D$  given by
\begin{eqnarray}
\lim_{\epsilon\searrow 0}\frac{1}{\omega_2\epsilon^2}\int_{\mathcal{H}_\epsilon(\bx)}\,|\bu(\by)-\bu(\bx)|\,d\by=0, 
\label{approx}
\end{eqnarray}
where $\mathcal{H}_\epsilon(\bx)$ is the ball of radius $\epsilon$ centered at $\bx$ and $\omega_2\epsilon^2$ is its area given in terms of the area of the unit disk $\omega_2$ times $\epsilon^2$.
The  jump set $\mathcal{J}_{\bu}$  for elements of  $SBD$ is defined to be the set of points of discontinuity which have two different one sided Lebesgue limits.  One sided Lebesgue limits of  $\bu$ with respect to a direction $\nu_{\bu(\bx)}$ are  denoted by $\bu^-(\bx)$, $\bu^+(\bx)$ and are given by
\begin{equation}
\begin{aligned}
&\lim_{\epsilon\searrow 0}\frac{1}{\epsilon^2\omega_2}\int_{\mathcal{H}^-_\epsilon(\bx)}\,|\bu(\by)-\bu^-(\bx)|\,d\by=0,\\
&\lim_{\epsilon\searrow 0}\frac{1}{\epsilon^2\omega_2}\int_{\mathcal{H}^+_\epsilon(\bx)}\,|\bu(\by)-\bu^+(\bx)|\,d\by=0,
\end{aligned}
\label{approxjump}
\end{equation}
where $\mathcal{H}^-_\epsilon(\bx)$ and $\mathcal{H}^+_\epsilon(\bx)$ are given by the intersection of $\mathcal{H}_\epsilon(\bx)$ with the half spaces $(\by-\bx)\cdot \bn_{\bu(\bx)}<0$ and $(\by-\bx)\cdot \bn_{\bu(\bx)}>0$ respectively. SBD functions have jump sets ${\mathcal{J}}_{\bu}$, described by a countable number of components $K_1,K_2,\ldots$, contained within smooth manifolds, with the exception of a set $K_0$ that has zero $1$ dimensional Hausdorff measure \\ 
\cite{AmbrosioCosicaDalmaso}. Here the notion of arc length is the one dimensional Hausdorff  measure of ${\mathcal J}_{\bu}$ and $\mathcal{H}^{d-1}({\mathcal J}_{\bu})=\sum_i\mathcal{H}^{1}(K_i)$.  
The  strain  of a displacement $\bu$ belonging to SBD, written as  $\mathcal{E}_{ij} \bu^0(t)=(\partial_{x_i}\bu^0_j+\partial_{x_j}\bu^0_i)/2$, is a generalization of the classic local strain tensor and is related to the 
nonlocal strain $S(\by,\bx,\bu^0)$  by 
\begin{equation}
\lim_{\epsilon\rightarrow 0}\frac{1}{\epsilon^2\omega_2}\int_{\mathcal{H}_\epsilon(x)}|{S(\by,\bx,\bu^0)}-\mathcal{E}\bu^0(\bx) \be\cdot \be|\,d\by=0,
\label{equatesandE}
\end{equation}
for almost every $\bx$ in $D$ with respect to $2$-dimensional Lebesgue measure $\mathcal{L}^2$. 
The symmetric part of the distributional derivative of $\bu$, $E \bu=1/2(\nabla \bu+\nabla \bu^T)$ for $SBD$ functions is a $2\times 2$ matrix valued Radon measure with absolutely continuous part described by the density $\mathcal{E}u$ and singular part described by the jump set \cite{AmbrosioCosicaDalmaso} and
\begin{eqnarray}
\langle E u,\Phi\rangle=\int_D\,\sum_{i,j=1}^d\mathcal{E}u_{ij}\Phi_{ij}\,d\bx+\int_{J_{u}}\,\sum_{i,j=1}^d(\bu^+_i - \bu^-_i)\bn_j\Phi_{ij}\,d\mathcal{H}^{1},
\label{distderiv}
\end{eqnarray}
for every continuous, symmetric matrix valued test function $\Phi$. In the sequel we will write $[\bu]=\bu^+-\bu^-$.
 Because $\bu^0$ has bounded Griffith energy \eqref{LEFMbound} we see that $\bu^0$ also belongs to $SBD^2$, that is the set of $SBD$ functions with square integrable strains $\mathcal{E}\bu$ and jump set with bounded Hausdorff $\mathcal{H}^1$ measure.


Now we establish Proposition \ref{measures}. To do so we first make the change of variables $\by=\bx+\epsilon{\xi}$ where ${\xi}$ belongs to the unit disk at the origin $\mathcal{H}_1(0)=\{|{\xi}|<1\}$ and $\be=\xi/|\xi|$. The strain is written
\begin{equation}
\begin{aligned}
\frac{\bu^\epsilon(\bx+\epsilon\xi)-\bu^\epsilon(\bx)}{\epsilon|\xi|}&:=D_{\be}^{\epsilon|\xi|}\bu^\epsilon,  \hbox{  and}\\
S(\by,\bx,\bu^\epsilon(t))&=D_{\be}^{\epsilon|\xi|}\bu^\epsilon\cdot\be, 
\end{aligned}
\label{straiin}
\end{equation}
and for infinitely differentiable scalar valued functions $\varphi$ and vector valued functions $\bw$ with compact support in $D$ we have
\begin{equation}\label{testgrad}
\lim_{\epsilon\rightarrow 0}D_{-\be}^{\epsilon|\xi|}\varphi=-\nabla\varphi\cdot\be,
\end{equation}
and
\begin{eqnarray}\label{teststrain}
\lim_{\epsilon\rightarrow 0}D_{\be}^{{\epsilon}|\xi|} \bw\cdot \be=\mathcal{E} \bw \,\be\cdot \be \label{grad}
\end{eqnarray}
where the convergence is uniform in $D$.
We now recall $S(\by,\bx,\bu^\epsilon(t))^-=D_{\be}^{\epsilon|\xi|}\bu^\epsilon\cdot\be^-$ defined by \eqref{decomposedetailsS}.
Here we point out that $D_{\be}^{\epsilon|\xi|}\bu^\epsilon\cdot\be^-=0$ for $\bx\in D$ and $\bx+\epsilon\xi\not\in D$. In this way $D_{\be}^{\epsilon|\xi|}\bu^\epsilon\cdot\be^-$ is well defined on $D\times\mathcal{H}_1(0)$.

As in inequality (6.73) of \cite{CMPer-Lipton} we have that
\begin{equation}\label{L2bound}
\begin{aligned}
\int_{D\times\mathcal{H}_1(0)}|\xi| J(|\xi|)|(D_{\be}^{{\epsilon}|\xi|}\bu^{\epsilon_k}\cdot \be)^-|^2\,d\xi\,d\bx<C,
\end{aligned}
\end{equation}
for all $\epsilon>0$.
From this we can conclude there exists a function $g(\bx,\xi)$ such that a subsequence $D_{\be}^{{\epsilon}|\xi|}\bu^{\epsilon}\cdot \be^-\rightharpoonup g(\bx,\xi)$ converges weakly in $L^2(D\times\mathcal{H}_1(\bx),\mathbb{R}^2)$ where the $L^2$ norm and inner product are with respect to the weighted measure $|\xi|J(|\xi|)d\xi d\bx $. Now for any positive number $\gamma$ and and any subset $D'$ compactly contained in $D$ we can argue as in (\cite{CMPer-Lipton} proof of Lemma 6.6) that $g(\bx,\xi)=\mathcal{E}\bu^0\be\cdot\be$ for all points in $D'$ with $|x_2|>\gamma$. Since $D'$ and $\gamma$ is arbitrary we get that
\begin{equation}\label{weaklim}
g(\bx,\xi)=\mathcal{E}\bu^0\be\cdot\be
\end{equation}
almost everywhere in  $D$. Additionally for any smooth scalar test function $\varphi(\bx)$ with compact support in $D$ straight forward computation gives 
\begin{equation}
\label{identify}
\begin{aligned}
\lim_{\epsilon\rightarrow0}\int_{D\times\mathcal{H}_1(0)}|\xi| J(|\xi|)D_{\be}^{{\epsilon}|\xi|}\bu^{\epsilon}\cdot \be^-\,d\xi\varphi(\bx)\,d\bx&=\int_{D\times\mathcal{H}_1(0)}|\xi| J(|\xi|)g(\bx,\xi)\,d\xi\varphi(\bx)\,d\bx\\
&=
\int_{D\times\mathcal{H}_1(0)}|\xi| J(|\xi|)\mathcal{E}\bu^0(\bx)\be\cdot\be\,d\xi\varphi(\bx)\,d\bx\\
&=C\int_D\,div\bu^0(\bx)\varphi(\bx)d\bx,
\end{aligned}
\end{equation}
Here $C=\omega_2\int_0^1r^2dr$ and we have used
\begin{equation}\label{convegDspoon}
\begin{aligned}
\frac{1}{\omega_2}\int_{\mathcal{H}_1(0)}|\xi|J(|\xi|)\be_ie_j\,d\xi=\delta_{ij}\int_0^1r^2J(r)\,dr.
\end{aligned}
\end{equation}
On the other hand for any smooth test function $\varphi$ with compact support in $D$ we can integrate by parts and use \eqref{testgrad} to write
\begin{equation}\label{L1}
\begin{aligned}
\lim_{\epsilon\rightarrow0}\int_{D\times\mathcal{H}_1(0)}|\xi| J(|\xi|)D_{\be}^{{\epsilon}|\xi|}\bu^{\epsilon}\cdot \be\varphi(\bx)\,d\xi\,d\bx&=\lim_{\epsilon\rightarrow0}\int_{D\times\mathcal{H}_1(0)}|\xi| J(|\xi|)D_{-\be}^{{\epsilon}|\xi|}\varphi(\bx)\bu^\epsilon\cdot \be,d\xi\,d\bx\\
&=-\int_{D\times\mathcal{H}_1(0)}|\xi| J(|\xi|)\bu^0\cdot \be\,\nabla\varphi(\bx)\cdot\be\,d\xi\,d\bx\\
&=-C\int_{D}\bu^0\cdot \nabla\varphi(\bx)\,d\bx\\
&=C\int_D\,tr{E\bu^0}\varphi(\bx)\,d\bx,
\end{aligned}
\end{equation}
where $E\bu^0$ is the strain of the $SBD^2$ limit displacement $\bu^0$.
Now since $\bu^0$ is in $SBD$ its weak derivitave satisfies \eqref{distderiv} and it follows on choosing $\Phi_{ij}={\delta}_{ij}\varphi$ that 
\begin{equation}\label{radon}
\begin{aligned}
\int_D\,tr{E\bu^0}\varphi\,d\bx=\int_D\,div\bu^0\varphi\,d\bx+\int_{\mathcal{J}_{\bu^0(t)}}[\bu^0]\cdot \bn\,\varphi\mathcal{H}^1(\bx),
\end{aligned}
\end{equation}
and note further  that
\begin{equation}\label{diff}
\begin{aligned}
\int_{D\times\mathcal{H}_1(0)}|\xi| J(|\xi|)D_{\be}^{{\epsilon}|\xi|}\bu^{\epsilon}\cdot \be\,d\xi\varphi(\bx)\,d\bx&=\int_{D\times\mathcal{H}_1(0)}|\xi| J(|\xi|)(D_{\be}^{{\epsilon}|\xi|}\bu^{\epsilon}\cdot \be)^-d\xi\varphi(\bx)\, d\bx\\
&+ \int_{D\times\mathcal{H}_1(0)}|\xi| J(|\xi|)(D_{\be}^{{\epsilon}|\xi|}\bu^{\epsilon}\cdot \be)^+d\xi \varphi(\bx)\,d\bx
\end{aligned}
\end{equation}
to conclude
\begin{equation}\label{radonlim}
\begin{aligned}
&\lim_{\epsilon \rightarrow 0}\int_{D\times\mathcal{H}_1(0)}|\xi| J(|\xi|)(D_{\be}^{{\epsilon}|\xi|}\bu^{\epsilon}\cdot \be)^+d\xi \varphi(\bx)\,d\bx\\
&=C\int_{\mathcal{J}_{\bu^0(t)}}[\bu^0]\cdot \bn\,\varphi\mathcal{H}^1(\bx).
\end{aligned}
\end{equation}
On changing variables we obtain the identities:
\begin{equation}\label{s1}
\begin{aligned}
&\lim_{\epsilon \rightarrow 0}\frac{1}{\epsilon^2}\int_{D}\int_{\mathcal{H}_\epsilon(\bx)}\frac{|\by-\bx|}{\epsilon} J^\epsilon(|\by-\bx|)
S(\by,\bx,\bu^\epsilon(t))^+\,d\by \,\varphi(\bx)\,d\bx\\
&=C\int_{\mathcal{J}_{\bu^0(t)}}[\bu^0]\cdot \bn\,\varphi\mathcal{H}^1(\bx).
\end{aligned}
\end{equation}
and
\begin{equation}
\label{s2}
\begin{aligned}
&\lim_{\epsilon\rightarrow0}\frac{1}{\epsilon^2}\int_{D}\int_{\mathcal{H}_\epsilon(\bx)}\frac{|\by-\bx|}{\epsilon} J^\epsilon(|\by-\bx|)
S(\by,\bx,\bu^\epsilon(t))^-\,d\by \,\varphi(\bx)\,d\bx\\
&=C\int_D\,div\bu^0(\bx)\varphi(\bx)d\bx.
\end{aligned}
\end{equation}
Since $S(\by,\bx,\bu^\epsilon(t))^+\not=0$ for $\bx$ in $SZ^\epsilon$ and zero otherwise  Proposition \ref{measures} is proved. We now use this proposition together with Hypotheses \ref{hyp3} to get Proposition \ref{cracksetjumpset}.

We fix $t$ and recall that the crack centerline is $C^\epsilon(t)=\{\ell(0)\leq x_1\leq \ell^\epsilon(t), \,x_2=0\}$. The sequence of numbers $\{\ell^\epsilon(t)\}_{\epsilon=1}^\infty$ is bounded so there exists at least one limit point and call it $\ell^0(t)$. In the sequel we will show that it is the only limit point of the sequence. 
Then there are at most two possibilities: a non-decreasing subsequence converging to $\ell^0$ or a non-increasing subsequence converging to $\ell^0$. We recall that  $SZ^\epsilon=F^\epsilon$ contains the thin rectangle $\{\ell(0)\leq x_1<\ell^\epsilon(t), \,-\epsilon<x_2<\epsilon\}$ and $\ell^\epsilon\rightarrow\ell^0$.
Then for either possibility if a positive smooth test function $\varphi$  with support set $supp\,\{\varphi\}\subset D$ intersects the interval $[0,\ell^0)$ on a set $B$ of nonzero one dimensional Lesbegue measure, i.e., $|B|>0$ then from Hypothesis \ref{hyp3} we have
\begin{equation}\label{hyp2playedout}
\begin{aligned}
&\lim_{\epsilon \rightarrow 0}\frac{1}{\epsilon^2\omega_2}\int_{SZ^\epsilon}\int_{\mathcal{H}_\epsilon(\bx)}\frac{|\by-\bx|}{\epsilon} J^\epsilon(|\by-\bx|)
S(\by,\bx,\bu^\epsilon(t))^+\,d\by \,\varphi(\bx)\,d\bx\\
&\geq K\alpha |B|>0,
\end{aligned}
\end{equation}
where $K$ is a constant depending on $\varphi$. Similarly if a positive smooth test function $\varphi$  with support set $supp\,\{\varphi\}\subset D$ intersects  $\mathcal{J}_{\bu^0(t)}$ on a set $B$ of nonzero one dimensional Lesbegue measure, i.e., $|B|>0$ then from Hypothesis \ref{hyp3} we also have
\begin{equation}\label{hyp2played}
\begin{aligned}
\int_{\mathcal{J}_{\bu^0(t)}}[\bu^0]\cdot \bn\,\varphi\mathcal{H}^1(\bx)\geq K\alpha |B|>0,
\end{aligned}
\end{equation}
where $K$ is a constant depending on $\varphi$.
We consider the non-decreasing case first. Suppose there is a positive smooth test function $\varphi$  with support set $supp\,\{\varphi\}\subset D$ intersecting $[0,\ell^0)$ on a nonzero set of one dimensional Lebesgue measure but not intersecting the jump set $\mathcal{J}_{\bu^0(t)}$. Then the left side of \eqref{convinmeasure} is positive but the right side is zero and there is a contradiction. Similarly we arrive at a contradiction if the support of a positive test function $\varphi$ intersects the jump set  $\mathcal{J}_{\bu^0(t)}$ on a set with positive one dimensional Lebesgue measure but not  $[0,\ell^0)$. So for this case we find that $|C^0\Delta \mathcal{J}_{\bu^0(t)}|=0$. One also easily arrives at contradictions for the non-increasing case as well. And we conclude that $\ell^0=|\mathcal{J}_{\bu^0(t)}|$. This also shows that all limit points coincide and Proposition \ref{cracksetjumpset} is proved.

\section{Convergence of nonlocal evolution to classic brittle fracture models}\label{s:proofofconvergence}

We start this section by establishing Proposition \ref{convergences}. The strong convergence 
\begin{equation}\label{convgences1}
\begin{aligned}
\bu^{\epsilon_n}\rightarrow \bu^0 & \hbox{  \rm strong in } C([0,T];\dot L^2(D;\mathbb{R}^2))
\end{aligned}
\end{equation}
follows immediately from the same arguments used to establish
Theorem 5.1 of \cite{CMPer-Lipton}.
The weak convergence
\begin{equation}\label{convgences2}
\begin{aligned}
\dot\bu^{\epsilon_n}\rightharpoonup \dot\bu^0 &\hbox{ \rm weakly in  }L^2(0,T;\dot L^2(D;\mathbb{R}^2))\end{aligned}
\end{equation}
follows noting that Theorem 2.2 of \cite{CMPer-Lipton} shows that
\begin{equation}
\begin{aligned}
\sup_{\epsilon_n>0}\int_0^T\Vert \dot\bu^{\epsilon_n}(t)\Vert^2_{L^2(D;\mathbb{R}^2)}dt<\infty.
\end{aligned}
\end{equation}
Thus we can pass to a subsequence also denoted by $\{\dot\bu^{\epsilon_n}\}_{n=1}^\infty$ that converges weakly to $\dot\bu^0$ in \\$L^2(0,T;\dot L^2(D;\mathbb{R}^2))$. To prove
\begin{equation}\label{convgences3}
\begin{aligned}
\ddot\bu^{\epsilon_n}\rightharpoonup \ddot\bu^0 &\hbox{ \rm weakly in }L^2(0,T;\dot H^1(D;\mathbb{R}^2)')
\end{aligned}
\end{equation}
we must show that
\begin{equation}\label{blf}
\sup_{\epsilon_n>0}\int_0^T\Vert\ddot\bu^{\epsilon_n}(t)\Vert^2_{\dot H^1(D;\mathbb{R}^2)'}\,dt<\infty,
\end{equation}
and existence of a weakly converging sequence follows.
To do this we consider the strong form of the evolution \eqref{energy based model2} which is an identity in $\dot L^2(D;\mathbb{R}^2)$ for all times $t$ in $[0,T]$. We multiply \eqref{energy based model2} with a test function $\bw$ from $\dot H^1(D;\mathbb{R}^2)$ and integrate over $D$. A straightforward integration  by parts gives 
\begin{equation}\label{energy based weakform est1}
\begin{aligned}
&\int_D\ddot{\bu}^{\epsilon_n}(\bx,t)\cdot {\bw}(\bx)d\bx\\
&=-\frac{1}{\rho}
\int_D \int_{H_{\epsilon_n}(\bx)\cap D} |\by-\bx|\partial_S\mathcal{W}^{\epsilon_n}(S(\by,\bx,\bu^{\epsilon_n}(t)))S(\by,\bx,\bw)\,d\by d\bx\\
&+\frac{1}{\rho}\int_D\bb(\bx,t)\cdot\bw(\bx)d\bx,
\end{aligned}
\end{equation}
and we now estimate the right hand side of \eqref{energy based weakform est1}.
For the first term on the righthand side we change variables $\by=\bx+\epsilon\xi$, $|\xi|<1$, with $d\by=\epsilon_n^2d\xi$  and write out $\partial_S\mathcal{W}^\epsilon(S(\by,\bx,\bu^\epsilon(t)))$ to get
\begin{equation}\label{est1}
\begin{aligned}
&I=-\frac{1}{\rho\omega_2}\int_{{D}\times{\mathcal{H}_1(0)}}\omega(\bx,\xi)|\xi|J(|\xi|)h'\left(\epsilon_n|\xi||D_{\be}^{\epsilon_n|\xi|}\bu^{\epsilon_n}\cdot \be|^2\right)\\
&\times2(D_{\be}^{\epsilon_n|\xi|}\bu^{\epsilon_n}\cdot \be) (D_{\be}^{\epsilon_n|\xi|}\bw\cdot \be)\,d\xi\,d\bx,
\end{aligned}
\end{equation}
where $\omega(\bx,\xi)$ is unity if $\bx+\epsilon\xi$ is in $D$ and zero otherwise. We define the sets
\begin{equation}\label{overstress}
\begin{aligned}
A^-_{\epsilon_n}&=\left\{(\bx,\xi) \hbox{ in }D\times\mathcal{H}_1(0);\,|D_{\be}^{\epsilon_n|\xi|}\bu^{\epsilon_n}\cdot\be|<\frac{{r}^c}{\sqrt{\epsilon_n|\xi|}}\right\}\\
A^+_{\epsilon_n}&=\left\{(\bx,\xi) \hbox{ in }D\times\mathcal{H}_1(0);\,|D_{\be}^{\epsilon_n|\xi|}\bu^{\epsilon_n}\cdot\be|\geq\frac{{r}^c}{\sqrt{\epsilon_n|\xi|}}\right\},
\end{aligned}
\end{equation}
with $D\times\mathcal{H}_1(0)=A^-_{\epsilon_n}\cup A^+_{\epsilon_n}$ and we write
\begin{equation}\label{i[plus}
I =I_1+I_2,
\end{equation}
where
\begin{equation}\label{est2}
\begin{aligned}
I_1=-\frac{1}{\rho\omega_2}\int_{D\times\mathcal{H}_1(0)\cap A_{\epsilon_n}^-}\omega(\bx,\xi)|\xi|J(|\xi|)h'\left(\epsilon_n|\xi||D_{\be}^{\epsilon_n|\xi|}\bu^{\epsilon_n}\cdot \be|^2\right)\\
\times 2(D_{\be}^{\epsilon_n|\xi|}\bu^{\epsilon_n}\cdot \be) (D_{\be}^{\epsilon_n|\xi|}\bw\cdot \be)\,d\xi\,d\bx,\\
I_2=-\frac{1}{\rho\omega_2}\int_{D\times\mathcal{H}_1(0)\cap A_{\epsilon_n}^+}\omega(\bx,\xi)|\xi|J(|\xi|)h'\left(\epsilon_n|\xi||D_{\be}^{\epsilon_n|\xi|}\bu^{\epsilon_n}\cdot \be|^2\right)\\
\times2(D_{\be}^{\epsilon_n|\xi|}\bu^{\epsilon_n}\cdot \be) (D_{\be}^{\epsilon_n|\xi|}\bw\cdot \be)\,d\xi\,d\bx,
\end{aligned}
\end{equation}
In what follows we introduce the the generic constant $C>0$ that is independent of $\bu^{\epsilon_n}$ and $\bw\in \dot H^1(D;\mathbb{R}^2)$. 
First note that $h$ is concave so $h'(r)$ is monotone decreasing for $r\geq 0$ and from Cauchy's inequality, and \eqref{L2bound} one has
\begin{equation}\label{est3}
\begin{aligned}
|I_1|&\leq\frac{2h'(0)C}{\rho\omega_2}\left(\int_{D\times\mathcal{H}_1(0)\cap A_{\epsilon_n}^-}|D_{\be}^{\epsilon_n|\xi|}\bw\cdot \be)|^2\,d\xi\,d\bx\right)^{1/2},\\
&\leq \frac{2h'(0)C}{\rho\omega_2}\left(\int_{\mathcal{H}_1(0)}\int_{D}|D_{\be}^{\epsilon_n|\xi|}\bw\cdot \be)|^2\,d\bx\,d\xi\right)^{1/2},
\end{aligned}
\end{equation}
For $t\in[0,T]$ the function $\bw$ is extended as an $H^1$ function to a larger domain $\tilde{D}$ containing $D$ such that there is a positive $\gamma$ such that  $0<\gamma<dist(D,\tilde{D})$ and $\Vert \bw\Vert_{H^1(\tilde{D};\mathbb{R}^2)}\leq C\Vert \bw\Vert_{H^1({D};\mathbb{R}^2)}$. For $\epsilon_n<\gamma$ the difference quotient satisfies
\begin{equation}\label{diffquo}
\Vert D_{\be}^{\epsilon_n|\xi|}\bw\cdot \be)\Vert_{L^2(D;\mathbb{R}^2)}\leq \Vert \bw\Vert_{H^1(\tilde{D};\mathbb{R}^2)}\leq C\Vert \bw\Vert_{H^1({D};\mathbb{R}^2)},
\end{equation}
for all $\xi\in\mathcal{H}_1(0)$ so
\begin{equation}\label{est4}
\begin{aligned}
|I_1|&\leq C\Vert \bw\Vert_{H^1({D};\mathbb{R}^2)}.
\end{aligned}
\end{equation}
Elementary calculation gives the estimate\\  (see equation (6.53) of \cite{CMPer-Lipton})
\begin{equation}\label{max}
\sup_{0\leq x<\infty}|h'({\epsilon_n}|\xi|{x}^2)2x|\leq\frac{2h'(\overline{r}^2){\overline{r}}}{\sqrt{\epsilon_n|\xi|}},
\end{equation}
and we also have (see equation (6.78) of \cite{CMPer-Lipton})
\begin{equation}
\int_{D\times\mathcal{H}_1(0)\cap A^{+}_{\epsilon_n}} \omega(\bx,\xi)J(|\xi|)\,d\xi\,d\bx<C\epsilon_n,
\label{Fprimesecond}
\end{equation}
so  Cauchy's inequality and the inequalities \eqref{diffquo}, \eqref{max}, \eqref{Fprimesecond}  give
\begin{equation}\label{est5}
\begin{aligned}
|I_2|&\leq\frac{1}{\rho\omega_2}\int_{D\times\mathcal{H}_1(0)\cap A_{\epsilon_n}^+}\omega(\bx,\xi)|\xi|J(|\xi|)\frac{2h'(\overline{r}^2){\overline{r}}}{\sqrt{\epsilon_n|\xi|}}|D_{\be}^{\epsilon_n|\xi|}\bw\cdot \be|\,d\xi\,d\bx,\\
&\leq\frac{1}{\rho\omega_2}\left(\int_{D\times\mathcal{H}_1(0)\cap A_{\epsilon_n}^+}\omega(\bx,\xi)|\xi|J(|\xi|)\frac{(2h'(\overline{r}^2){\overline{r}})^2}{{\epsilon_n|\xi|}}\,d\xi\,d\bx\right)^{1/2}\times\\
&\left(\int_{D\times\mathcal{H}_1(0)\cap A_{\epsilon_n}^+}\omega(\bx,\xi)|\xi|J(|\xi|)|D_{\be}^{\epsilon_n|\xi|}\bw\cdot \be|^2\,d\xi\,d\bx\,dt\right)^{1/2}\\
&\leq C\Vert \bw\Vert_{H^1({D};\mathbb{R}^2)},
\end{aligned}
\end{equation}
and we conclude that the first term on the right hand side of \eqref{energy based weakform est1} admits the estimate
\begin{equation}\label{fin}
\begin{aligned}
|I| &\leq |I_1|+|I_2| \leq  C\Vert \bw\Vert_{H^1({D};\mathbb{R}^2)},
\end{aligned}
\end{equation}
for all $\bw\in H^1(D;\mathbb{R}^2)$.

It follows immediately from \eqref{l2} that the second term on the right hand side of  \eqref{energy based weakform est1} satisfies the estimate
\begin{equation}\label{blf2}
\begin{aligned}
\frac{1}{\rho}\left|\int_D\bb(\bx,t)\cdot\bw(\bx)\,d\bx\right |\leq C\Vert\bw\Vert_{H^1(D;\mathbb{R}^2)},\hbox{ for all $\bw\in H^1(D;\mathbb{R}^2)$}
\end{aligned}
\end{equation}
From \eqref{fin} and \eqref{blf2} we conclude that there exists a $C>0$ so that
\begin{eqnarray}\label{energy based weakform est3}
\begin{aligned}
\left |\int_D\ddot{\bu}^{\epsilon_n}(\bx,t)\cdot {\bw}(\bx)\,d\bx\right |\leq C\Vert\bw\Vert_{H^1(D;\mathbb{R}^2)},\hbox{ for all $\bw\in \dot H^1(D;\mathbb{R}^2)$}
\end{aligned}
\end{eqnarray}
so
\begin{equation}\label{bounded}
\begin{aligned}
&\sup_{\epsilon_n>0}\sup_{t\in [0,T]}\frac{\int_D\ddot{\bu}^{\epsilon_n}(\bx,t)\cdot {\bw}(\bx)d\bx}{\Vert\bw\Vert_{H^1(D;\mathbb{R}^2)}}<C,\hbox{ for all $\bw\in \dot H^1(D;\mathbb{R}^2)$},
\end{aligned}
\end{equation}
and \eqref{blf} follows. The estimate \eqref{blf} implies weak compactness and passing to subsequences if necessary we deduce that
$\ddot\bu^{\epsilon_n}\rightharpoonup \ddot\bu^0 $ weakly in $L^2(0,T;\dot H^1(D;\mathbb{R}^2)')$ and Proposition \ref{convergences} is proved.

To establish Proposition \ref{momentumlim} we first note that it is easy to see that $\ddot\bu^0$ also belongs to  $L^2(0,T; \,H^{1,0}(D\setminus C^0(t),\mathbb{R}^2)')$ from Proposition  \ref{convergences} and we show that $\bu^0$ is a solution of \eqref{momentumlimit}.
We take a test function $\bw(\bx)$ that is infinitely differentiable on $D$ with support set that does not intersect the crack fix $t\in[0,T]$. Multiplying \eqref{energy based model2} by this test function and integration by parts gives as before
\begin{equation}\label{energy based weakform est15}
\begin{aligned}
&\rho\int_D\ddot{\bu}^{\epsilon_n}(\bx,t)\cdot {\bw}(\bx)d\bx\\
 &=-\int_D \int_{\mathcal{H}_{\epsilon_n}(\bx)\cap D} |\by-\bx|\partial_S\mathcal{W}^{\epsilon_n}(S(\by,\bx,\bu^{\epsilon_n}(t)))S(\by,\bx,\bw)\,d\by d\bx\\
&+\int_D\bb(\bx,t)\cdot\bw(\bx)d\bx,
\end{aligned}
\end{equation}
The goal is to pass to the $\epsilon_n=0$  in this equation to recover \eqref{momentumlimit}.
Using arguments identical to those above we find that for fixed $t$ that on passage to a possible subsequence also denoted by $\{\epsilon_n\}$  one recovers the term on the left hand side of \eqref{momentumlimit}, i.e.,
\begin{equation}\label{momentumlimitfirstterm}
\begin{aligned}
\lim_{\epsilon_n\rightarrow 0}\rho\int_D\ddot{\bu}^{\epsilon_n}(\bx,t)\cdot {\bw}(\bx)d\bx &=\rho\langle \ddot\bu^0,\bw\rangle.
\end{aligned}
\end{equation}
To recover the $\epsilon_n=0$ limit of the first term on the right hand side of  \eqref{energy based weakform est15} we see that \eqref{teststrain} and \eqref{weaklim}  hold and
identical arguments as in the proof of Theorem 6.7 of \cite{CMPer-Lipton3} show that on passage to a further subsequence if necessary one obtains
\begin{equation}\label{limitmomentumbalance 1st termrightside}
\begin{aligned}
&-\lim_{\epsilon_n\rightarrow 0}\int_D \int_{H_{\epsilon_n}(\bx)\cap D} |\by-\bx|\partial_S\mathcal{W}^{\epsilon_n}(S(\by,\bx,\bu^{\epsilon_n}(t)))S(\by,\bx,\bw)\,d\by d\bx\\
&=-\int_{D}\mathbb{C}\mathcal{E}\bu^0:\mathcal{E}\bw\,\,d\bx,\\
&=-\int_{D\setminus C^0(t)}\mathbb{C}\mathcal{E}\bu^0:\mathcal{E}\bw\,\,d\bx,
\end{aligned}
\end{equation}
where the last equality follows since the support of $\bw$ is away from the crack. 
The second term on the right hand side of  \eqref{energy based weakform est15} is a bounded linear functional on $H^{1,0}(D\setminus C^0(t),\mathbb{R}^2)$ and we make the identification $\int_D\bb\cdot\bw\,dx=\langle \bb,\bw\rangle$ and the last term on the righthand side of \eqref{momentumlimit} follows. 
This shows that \eqref{momentumlimit} holds for all infinitely differentiable test functions with support away from the crack and Proposition \ref{momentumlim} now follows from density of the test functions in $H^{1}(D\setminus C^0(t),\mathbb{R}^2)$.

To establish Proposition \ref{traction2} we first show that $\ddot\bu^0(t)$ is a bounded linear functional on the spaces  \,$H^{1,0}(Q^{\pm}_\beta(t),\mathbb{R}^2)$ for a.e. $t\in [0,T]$. We recall $\ell^{\epsilon_n}(t)\rightarrow\ell^0(t)$ and for $\beta$ such that $\ell(0)\leq\ell^0(t)-\beta$ we only consider $\epsilon_n$ so that $\ell^{0}(t)-\beta/2<\ell^{\epsilon_n}(t)$ and $\epsilon_n<\beta/2$. We make this choice since the interval $\{\ell(0)\leq x_1<\ell^0(t)-\beta;\, x_2=0\}$ is now included in the crack center line $C^{\epsilon_n}(t)$ see Definition \ref{def1}.
We multiply \eqref{energy based model2} with a test function $\bw$ from $H^{1,0}(Q^{\pm}_\beta(t),\mathbb{R}^2)$ and integrate over $D$ and perform a straight forward integration  by parts to get
\begin{equation}\label{energy based weakform est 20}
\begin{aligned}
&\int_{Q^{\pm}_\beta(t)}\ddot{\bu}^{\epsilon_n}(\bx,t)\cdot {\bw}(\bx)d\bx\\
 &=-\frac{1}{\rho}
\int_{Q^{\pm}_\beta(t)} \int_{H_{\epsilon_n}(\bx)\cap Q^{\pm}_\beta(t)} |\by-\bx|\partial_S\mathcal{W}^{\epsilon_n}(S(\by,\bx,\bu^{\epsilon_n}(t)))S(\by,\bx,\bw)\,d\by d\bx
\end{aligned}
\end{equation}
We can bound the term on the righthand side of \eqref{energy based weakform est 20} using the same arguments used to bound \eqref{est1}. The only difference is in the extension of the test function from the rectangles $Q^+_\beta(t)$ or $Q_\beta^-(t)$ to larger rectangles. Here, given a fixed  $\gamma>0$ for $t\in[0,T]$ the function $\bw$ is extended as an $H^1$ function on $Q^{+}_\beta(t)$ to the larger rectangle $\tilde{Q}^{+}_\beta(t)$ containing $Q^{+}_\beta(t)$ given by  $\tilde{Q}^{+}_\beta(t)=\{\ell(0)<x_1<\l^0(t)-\beta;-\gamma<x_2<b/2-\delta\}$ and $\Vert \bw\Vert_{H^1(\tilde{Q}^{+}_\beta(t),\mathbb{R}^2)}\leq C\Vert \bw\Vert_{H^1(Q^{+}_\beta(t),\mathbb{R}^2)}$.  Similarly for $t\in[0,T]$ the function $\bw$ is extended as an $H^1$ function on $Q^{-}_\beta(t)$ to the larger rectangle $\tilde{Q}^{-}_\beta(t)$ containing $Q^{-}_\beta(t)$ given by  $\tilde{Q}^{-}_\beta(t)=\{\ell(0)<x_1<\l^0(t)-\beta;-b/2+\delta<x_2<\gamma\}$ and $\Vert \bw\Vert_{H^1(\tilde{Q}^{-}_\beta(t),\mathbb{R}^2)}\leq C\Vert \bw\Vert_{H^1(Q^{-}_\beta(t),\mathbb{R}^2)}$. For our choice of $\epsilon_n$ the difference quotient satisfies
\begin{equation}\label{diffquo2}
\Vert D_{\be}^{\epsilon_n|\xi|}\bw\cdot \be)\Vert_{L^2(Q^{\pm}_\beta,\mathbb{R}^2)}\leq \Vert \bw\Vert_{H^1(\tilde{Q}^{\pm}_\beta,\mathbb{R}^2)}\leq C\Vert \bw\Vert_{H^1(Q^{\pm}_\beta,\mathbb{R}^2)},
\end{equation}
for all $\xi\in\mathcal{H}_1(0)$. We can then proceed as before to find for a.e. $t\in [0,T]$ that
\begin{equation}\label{boundedtfixed}
\begin{aligned}
 \sup_{\epsilon_n>0}\Vert\ddot\bu^{\epsilon_n}(t)\Vert^2_{H^1(Q^{\pm}_\beta(t),\mathbb{R}^2)'}<\infty,
\end{aligned}
\end{equation}
The estimate \eqref{boundedtfixed}
 implies compactness with respect to weak convergence and passing to subsequences if necessary we deduce that
$\ddot\bu^{\epsilon_n}(t)\rightharpoonup \ddot\bu^0(t)$ weakly in $H^1(Q^{\pm}_\beta(t),\mathbb{R}^2)'$ and we see that $\ddot\bu^0(t)$ 
 is a bounded linear functional on the spaces $H^1(Q^{\pm}_\beta(t),\mathbb{R}^2)$ for a.e. $t\in [0,T]$.  
 
To illustrate the ideas we now recover \eqref{momentumlimit2} on $Q^+_\beta(t)$. We first consider \eqref{energy based weakform est 20} with infinitely differentiable test functions $\bw(\bx)$  on $Q^+_\beta(t)$ with support sets that do not intersect the sets $\{\ell(0)\leq x_1<\ell(t)-\beta\,;,x_2=b/2-\delta\}$, $\{x_1=\ell^0(t)-\beta;\, 0<x_2<b/2-\delta\}$, and $\{x_1=\ell(0);\, 0<x_2<b/2-\delta\}$.
Passing to subsequences as necessary we recover the limit equation \eqref{momentumlimit2} using the same arguments that were used to pass to the limit in \eqref{momentumlimit}.  Proposition \ref{traction2} now follows using the density of these trial fields in $H^1(Q^+_\beta,\mathbb{R}^2)$. An identical procedure works for the rectangles $Q^-_\beta$ and Proposition \ref{traction2} is proved.

\section{Power balance on subdomains containing the crack tip for the nonlocal model}\label{s:energyrate}

\begin{figure} 
\centering
\begin{tikzpicture}[xscale=0.60,yscale=0.60]

\draw [-,thick] (-2,0) -- (-0.0,0);

\draw [thick] (-2,-3) rectangle (2,3);

\draw [thick] (-1,-1) rectangle (1,1);

\node [above] at (0.0,0.0) {$P^{\epsilon}_\delta$};

\node [right] at (.9,0.0) {$\Gamma^{\epsilon}_\delta$};

\draw [->,thick] (.5,-1.5) -- (1.5,-1.5);

\node [left] at (0.5,-1.5) {$V^{\epsilon}$};

\end{tikzpicture} 
\caption{{ \bf Contour $\Gamma^{\epsilon}_\delta$ surrounding the domain $P^{\epsilon}_\delta$ moving with the crack centerline tip velocity $V^{\epsilon}$.}}
 \label{Pepsilon}
\end{figure}
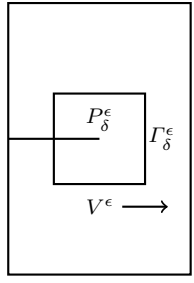

In this section we derive the power balance for the nonlocal model using the balance of linear momentum given by \eqref{energy based model2}. Consider the rectangular contour $\Gamma^\epsilon_\delta(t)$ of diameter $\delta$ surrounding the domain $\mathcal{P}^\epsilon_\delta(t)$ containing the crack tip. We suppose $\mathcal{P}^\epsilon_\delta(t)$ is moving with the crack tip speed $V^\epsilon(t)=\dot\ell^\epsilon(t)$ see Figure \ref{Pepsilon}. It will be shown that the rate of change of energy  inside  $\mathcal{P}^\epsilon_\delta(t)$  for the nonlocal dynamics is given by \eqref{powerbfornonloc}.
We start by introducing a nonlocal divergence theorem applied to the case at hand. To expedite taking $\epsilon\rightarrow0$ limits we make the change of variables $\by=\bx+\epsilon\xi$ where $\xi\in\mathcal{H}_1(0)$. The strain $S(\by,\bx,\bu^\epsilon(t))$ transforms to $D_{\be}^{\epsilon|\xi|}\bu^\epsilon\cdot\be$ and the work done in straining the material between points $\by$ and $\bx$ given by $|\by-\bx|\partial_S\mathcal{W}^\epsilon(S(\by,\bx,\bu^\epsilon(t)))$ transforms in the new variables to
\begin{equation}\label{newvlbls}
\epsilon|\xi|\partial_S\mathcal{W}^\epsilon(D_{\be}^{\epsilon|\xi|}\bu^\epsilon\cdot\be)=\frac{2|\xi|J(|\xi|)}{\epsilon^2\omega_2}h'(\epsilon|\xi||D_{\be}^{\epsilon|\xi|}\bu^\epsilon\cdot\be|^2)D_{\be}^{\epsilon|\xi|}\bu^\epsilon\cdot\be.
\end{equation}
The nonlocal divergence theorem is given by
\begin{proposition}\label{divg}
\begin{equation}\label{divthm}
\begin{aligned}
&\epsilon^2\int_{P^{\epsilon}_\delta(t)}\int_{\mathcal{H}_1(0)}D_{-\be}^{{\epsilon}|\xi|}\left[ \epsilon|\xi|\partial_S\mathcal{W}^\epsilon(D_{\be}^{{\epsilon}|\xi|}\bu^\epsilon\cdot\be)\bw\cdot\be\right]\,d\xi d\bx\\
&=\epsilon^2\int_{H_1(0)}\int_{(P^{\epsilon}_\delta(t)-\epsilon\xi)\setminus P^{\epsilon}_\delta(t)}\partial_S\mathcal{W}^\epsilon(D_{\be}^{{\epsilon}|\xi|}\bu^\epsilon\cdot\be)\bw\cdot\be\,d\bx d\xi\\
&-\epsilon^2\int_{H_1(0)}\int_{P^{\epsilon}_\delta(t)\setminus (P^{\epsilon}_\delta(t)-\epsilon\xi)}\partial_S\mathcal{W}^\epsilon(D_{\be}^{{\epsilon}|\xi|}\bu^\epsilon\cdot\be)\bw\cdot\be\,d\bx d\xi.
\end{aligned}
\end{equation}
\end{proposition}
This identity follows on applying the definition of $D_{-\be}^{\epsilon|\xi|}\varphi=(\varphi(\bx-\epsilon\xi)-\varphi(\bx))/\epsilon|\xi|$ for scalar fields  $\varphi$ and Fubini's theorem. When convenient we set $A^{\epsilon}_\delta(t)=D\setminus P^\epsilon_\delta(t)$ and rewrite the last two terms of \eqref{divthm} in  $\bx$ and $\by$ variables to get
\begin{equation}\label{divthm2}
\begin{aligned}
\epsilon^2&\int_{P^{\epsilon}_\delta(t)}\int_{\mathcal{H}_1(0)}D_{-\be}^{{\epsilon}|\xi|}\left[ \epsilon|\xi|\partial_S\mathcal{W}^\epsilon(D_{\be}^{{\epsilon}|\xi|}\bu^\epsilon\cdot\be)\bw\cdot\be\right]\,d\xi d\bx\\
&=\int_{A^{\epsilon}_\delta(t)}\int_{\mathcal{H}_{\epsilon}(\bx)\cap P^{\epsilon}_\delta(t) }\partial_S\mathcal{W}^{\epsilon}(S(\by,\bx,\bu^\epsilon(t)))(\bw(\bx)+\bw(\by))\cdot\be_{\by-\bx}\,d\by d\bx.
\end{aligned}
\end{equation}
Finally to expedite a convenient parameterization for passing to the $\epsilon_n=0$  limit we rewrite \eqref{divthm2} in $\bx$ and $\xi$ variables to get
\begin{equation}\label{divthm3}
\begin{aligned}
&\epsilon^2\int_{P^{\epsilon}_\delta(t)}\int_{\mathcal{H}_1(0)}D_{-\be}^{{\epsilon}|\xi|}\left[ \epsilon|\xi|\partial_S\mathcal{W}^\epsilon(D_{\be}^{{\epsilon}|\xi|}\bu^\epsilon\cdot\be)\bw\cdot\be\right]\,d\xi d\bx\\
&=\epsilon^2\int_{\mathcal{H}_{1}(0)}\int_{(P^{\epsilon}_\delta(t)-\epsilon_n\xi)\setminus P^{\epsilon}_\delta(t) }\partial_S\mathcal{W}^{\epsilon}(D_{\be}^{{\epsilon}|\xi|}\bu^\epsilon\cdot\be)(\bw(\bx)+\bw(\bx+\epsilon_n\xi))\cdot\be\,d\bx d\xi.
\end{aligned}
\end{equation}

We have the product rule given by
\begin{equation}\label{product}
\begin{aligned}
&\epsilon^2\int_{P^{\epsilon}_\delta(t)}\int_{\mathcal{H}_1(0)}D_{-\be}^{{\epsilon}|\xi|}\left[ \epsilon|\xi|\partial_S\mathcal{W}^\epsilon(D_{\be}^{{\epsilon}|\xi|}\bu^\epsilon\cdot\be)\bw\cdot\be\right]\,d\xi d\bx\\
&=-\epsilon^2\int_{P^{\epsilon}_\delta(t)}\int_{\mathcal{H}_1(0)}2\partial_S\mathcal{W}^\epsilon(D_{\be}^{{\epsilon}|\xi|}\bu^\epsilon\cdot\be)\be\cdot\bw\,d\xi d\bx\\
&-\epsilon^2\int_{P^{\epsilon}_\delta(t)}\int_{\mathcal{H}_1(0)}\epsilon|\xi|\partial_S\mathcal{W}^\epsilon(D_{\be}^{{\epsilon}|\xi|}\bu^\epsilon\cdot\be)D_{\be}^{{\epsilon}|\xi|}\bw\cdot\be\,d\xi d\bx.
\end{aligned}
\end{equation}

We now derive Proposition \ref{changinternal}.
Multiplying both sides of \eqref{energy based model2} by $\dot\bu^{\epsilon_n}$, integration over $P^{\epsilon_n}_\delta(t)$, and applying the product rule gives
\begin{equation}\label{balance2}
\begin{aligned}
&\int_{P^{\epsilon_n}_\delta(t)}\partial_t\frac{\rho|{\bu}^{\epsilon_n}|^2}{2}\,d\bx\\
&=\epsilon_n^2\int_{P^{\epsilon_n}_\delta(t)}\int_{\mathcal{H}_1(0)}2\partial_S\mathcal{W}^{\epsilon_n}(D_{\be}^{{\epsilon_n}|\xi|}\bu^\epsilon\cdot\be)\dot\bu^{\epsilon_n}\cdot\be\,d\xi d\bx\\
&=-\epsilon_n^2\int_{P^{\epsilon_n}_\delta(t)}\int_{\mathcal{H}_1(0)}D_{-\be}^{{\epsilon_n}|\xi|}\left[ \epsilon|\xi|\partial_S\mathcal{W}^{\epsilon_n}(D_{\be}^{{\epsilon_n}|\xi|}\bu^{\epsilon_n}\cdot\be)\dot\bu^{\epsilon_n}\cdot\be\right]\,d\xi d\bx\\
&-\epsilon_n^2\int_{P^{\epsilon_n}_\delta(t)}\int_{\mathcal{H}_1(0)}\epsilon_n|\xi|\partial_S\mathcal{W}^{\epsilon_n}(D_{\be}^{{\epsilon_n}|\xi|}\bu^{\epsilon_n}\cdot\be)D_{\be}^{{\epsilon_n}|\xi|}\dot\bu^{\epsilon_n}\cdot\be\, d\xi d\bx
\end{aligned}
\end{equation}
Define the energy density
\begin{equation}\label{density}
W^{\epsilon_n}(\bx,t)=\epsilon_n^2\int_{\mathcal{H}_1(0)}\epsilon_n|\xi|\mathcal{W}^{\epsilon_n}(D_{\be}^{{\epsilon_n}|\xi|}\bu^{\epsilon_n}\cdot\be)\,d\xi.
\end{equation}
We observe that the change in energy density with respect to time is given by
\begin{equation}\label{balancerate}
\begin{aligned}
\dot W^{\epsilon_n}&=\int_{\mathcal{H}_1(0)}\epsilon_n^3|\xi|\partial_S\mathcal{W}^{\epsilon_n}(D_{\be}^{{\epsilon_n}|\xi|}\bu^{\epsilon_n}\cdot\be)D_{\be}^{{\epsilon_n}|\xi|}\dot\bu^{\epsilon_n}\cdot\be\, d\xi,
\end{aligned}
\end{equation}
and \eqref{balance2} becomes
\begin{equation}\label{balance3}
\begin{aligned}
&\int_{P^{\epsilon_n}_\delta(t)}\,\dot T^{\epsilon_n}+\dot W^{\epsilon_n}\,d\bx\\
&=-\int_{P^{\epsilon_n}_\delta(t)}\int_{\mathcal{H}_1(0)}\epsilon_n^2D_{-\be}^{{\epsilon_n}|\xi|}\left[ \epsilon|\xi|\partial_S\mathcal{W}^{\epsilon_n}(D_{\be}^{{\epsilon_n}|\xi|}\bu^{\epsilon_n}\cdot\be)\dot\bu^{\epsilon_n}\cdot\be\right]\,d\xi d\bx,
\end{aligned}
\end{equation}
where $\dot T^{\epsilon}=\partial_t(\rho|\dot\bu^{\epsilon_n}|^2/2)$. 

The domain $P^{\epsilon_n}_\delta(t)$ is traveling with the crack velocity $V^{\epsilon_n}(t)$ so Reynolds transport theorem together with Proposition \ref{divg} and \eqref{divthm2} deliver the power balance:
 \begin{equation}\label{balance4}
\begin{aligned}
&\frac{d}{dt}\int_{P^{\epsilon_n}_\delta(t)}\,T^{\epsilon_n}+ W^{\epsilon_n}\,d\bx\\
&=\int_{\partial P^{\epsilon_n}_\delta(t)}\,(T^{\epsilon_n}+ W^{\epsilon_n})V^{\epsilon_n}\be^1\cdot\bn\,\,ds\\
&-\int_{A^{\epsilon}_\delta(t)}\int_{\mathcal{H}_{\epsilon}(\bx)\cap P^{\epsilon}_\delta(t) }\partial_S\mathcal{W}^{\epsilon}(S(\by,\bx,\bu^\epsilon(t)))\,\be_{\by-\bx}\cdot(\dot\bu^{\epsilon_n}(\bx)+\dot\bu^{\epsilon_n}(\by))\,d\by d\bx
\end{aligned}
\end{equation}
and \eqref{powerbfornonloc} follows and Proposition \ref{changinternal} is proved.

\section{Crack tip motion and power balance for the local model}\label{s:Jintegral}
In this section we complete the proof of Proposition \ref{limitspowerflow} and establish \eqref{nonloctolocrate1}. We start by establishing the crucial identity
\begin{equation}\label{nonlocSurfaceterm}
\begin{aligned}
 \int_{\Gamma^{\epsilon_n}_\delta(t)}\,(T^{\epsilon_n}+W^{\epsilon_n})V^{\epsilon_n}\be^1\cdot\bn\,ds&=-\mathcal{G}_c V^{\epsilon_n}(t)+O(\delta). 
 \end{aligned}
\end{equation}

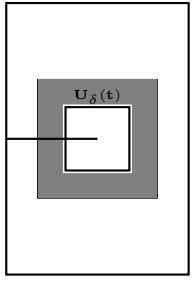
\begin{figure} 
\centering
\begin{tikzpicture}[xscale=0.60,yscale=0.60]


\draw [thick] (-2,-3) rectangle (2,3);

\draw [thick] (-1.3,-1.3) rectangle (1.3,1.3);

\begin{scope}
\clip [postaction={draw=gray, line width=6.5mm}] (-1.3,-1.3) rectangle + (2.6,2.6);
\end{scope}
\draw [thick] (-0.7,-0.7) rectangle (0.7,0.7);

\node [above] at (0.0,0.6) {$\scriptscriptstyle{\bf{U_\delta(t)}}$};

\draw [-,thick] (-2,0) -- (-0.0,0);




\end{tikzpicture} 
\caption{{ $U_\delta(t)$.}}
 \label{oshape}
\end{figure}

Since $V^{\epsilon_n}\be^1\rightarrow V\be^1$ we note that the contours $\Gamma^{\epsilon_n}_\delta$  converge to $\Gamma_\delta(t)$ and for $\epsilon_n$ small enough that $\Gamma^{\epsilon_n}_\delta$ lie in a  ``O'' shaped domain $U_\delta(t)$ surrounding the crack tip with boundary  at least $\delta/2$ away from the tip see figure \ref{oshape}.  So we can conclude from the hypothesis of Proposition \ref{limitspowerflow} that $\bu^{\epsilon_n}(t)$, $\dot\bu^{\epsilon_n}(t)$, $S(\by,\bx,\bu^{\epsilon_n})$  converge uniformly to $\bu^0(t)$,  $\dot\bu^0$(t), and $\mathcal{E}\bu^0\be\cdot\be$ on $U_\delta(t)$ for $t\in[0,T]$.

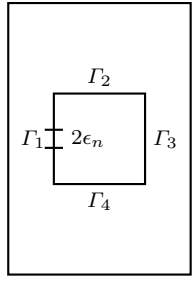
\begin{figure} 
\centering
\begin{tikzpicture}[xscale=0.60,yscale=0.60]

\draw [thick] (-2,-3) rectangle (2,3);

\draw [thick] (-1,-1) rectangle (1,1);

\draw [thick] (-1.2,.2) -- (-.8,.2);

\draw [thick] (-1.2,-.2) -- (-.8,-.2);

\node [right] at (-0.8,0.0) {$2\epsilon_n$};

\node [above] at (0,1) {$\Gamma_2$};

\node [right] at (1,0.0) {$\Gamma_3$};

\node [left] at (-1,0) {$\Gamma_1$};

\node [below] at (0,-1) {$\Gamma_4$};

\end{tikzpicture} 
\caption{{ \bf Contour $\Gamma^{\epsilon_n}_\delta$ split into four sides and $\Gamma_1$ divided into one part away from crack centerline and two parts close to the crack centerline.}}
 \label{sidess}
\end{figure}

Now we denote the four sides of the rectangular contour $\Gamma^{\epsilon_n}_\delta$ by $\Gamma_i$, $i=1,\ldots,4$ in figure \ref{sidess}. There is no contribution of the integrand to the integral on the lefthand side of \eqref{nonlocSurfaceterm} on sides $2$ and $4$ as $\be^1\cdot\bn=0$ there. On side $3$ it follows from the uniform convergence on $U_\delta(t)$ and \eqref{LEFMequality} that both  $|T^{\epsilon_n}|<C$ and $|W^{\epsilon_n}|<C$  on $U_\delta(t)$ for all sufficiently small $\epsilon_n$ so 
\begin{equation}\label{gamma3}
\left\vert \int_{\Gamma_3}\,(T^{\epsilon_n}+W^{\epsilon_n})V^{\epsilon_n}\be^1\cdot\bn\,ds\right\vert=O(\delta).
\end{equation}
On side $1$ we partition the contour $\Gamma_1$ into three parts. The first part is given by all points on $\Gamma_1$ that are further than $\epsilon_n$ away from $x_2=0$ call this $\Gamma_{1,+}$ and as before
\begin{equation}\label{gamma1plus}
\left\vert \int_{\Gamma_{1,+}}\,(T^{\epsilon_n}+W^{\epsilon_n})V^{\epsilon_n}\be^1\cdot\bn\,ds\right\vert=O(\delta).
\end{equation}
The part of $\Gamma_1$ with $0\leq x_2\leq\epsilon_n$ is denoted $\Gamma_1^+$ and the part with $-\epsilon_n\leq x_2<0$ is denoted $\Gamma_1^-$. Now we calculate
\begin{eqnarray}
\begin{aligned}
&\int_{\Gamma_1^+}W^{\epsilon_n}V^{\epsilon_n}\be^1\cdot\bn\,ds\\
&=-V^{\epsilon_n}\int_{\Gamma_1^+}W^{\epsilon_n}\,ds\\
&=-V^{\epsilon_n}\int_{\Gamma_1^+}\int_{\mathcal{H}_\epsilon(\bx)\cap K_{\epsilon_n}^+}|\by-\bx|\mathcal{W}^\epsilon (S(\by,\bx,\bu^{\epsilon_n}(t)))\,d\by ds\\
&-V^{\epsilon_n}\int_{\Gamma_1^+}\int_{\mathcal{H}_\epsilon(\bx)\cap K_{\epsilon_n}^-}|\by-\bx|\mathcal{W}^\epsilon (S(\by,\bx,\bu^{\epsilon_n}(t)))\,d\by ds
\end{aligned}
\label{formulaplus}
\end{eqnarray}
Here $\mathcal{H}_\epsilon(\bx)\cap K_{\epsilon_n}^+$ is the subset of $\by$ in $\mathcal{H}_{\epsilon_n}(\bx)$ for which the vector with end points $\by$ and $\bx$ crosses the crack centerline and $\mathcal{H}_\epsilon(\bx)\cap K_{\epsilon_n}^-$ is the subset of $\by$ in $\mathcal{H}_{\epsilon_n}(\bx)$ for which  the vector with end points $\by$ and $\bx$ does not cross the crack centerline. From Hypothesis \ref{hyp2} and calculation as in section \ref{toughelastic} gives
\begin{eqnarray}
\begin{aligned}
&\int_{\Gamma_1^+}\int_{\mathcal{H}_\epsilon(\bx)\cap K_{\epsilon_n}^+}|\by-\bx|\mathcal{W}^\epsilon (S(\by,\bx,\bu^{\epsilon_n}(t)))\,d\by ds\\
&=2\int_0^\epsilon\int_z^\epsilon\int_0^{\arccos(z/\zeta)}\mathcal{W}^\epsilon(\mathcal{S}_+)\zeta^2\,d\psi\,d\zeta\,dz\\
&=\frac{\mathcal{G}_c}{2},
\end{aligned}
\label{formulaplusmoreG}
\end{eqnarray}
and it follows from the  uniform convergence on $U_\delta(t)$ and calculating as in \eqref{LEFMequality} we get that
\begin{eqnarray}
\begin{aligned}
\left\vert\int_{\Gamma_1^+}\int_{\mathcal{H}_\epsilon(\bx)\cap K_{\epsilon_n}^-}|\by-\bx|\mathcal{W}^\epsilon (S(\by,\bx,\bu^{\epsilon_n}(t)))\,d\by ds\right\vert=&O(\delta).
\end{aligned}
\label{formulaplusmoreGd}
\end{eqnarray}
From \eqref{formulaplusmoreG} and \eqref{formulaplusmoreGd} we conclude that
\begin{eqnarray}
\begin{aligned}
\int_{\Gamma_1^+}W^{\epsilon_n}V^{\epsilon_n}\be^1\cdot\bn\,ds=&-V^{\epsilon_n}\int_{\Gamma_1^+}W^{\epsilon_n}\,ds\\
=-V^{\epsilon_n}\frac{\mathcal{G}_c}{2}+O(\delta).
\end{aligned}
\label{formulaplusplus}
\end{eqnarray}
An identical calculation shows 
\begin{eqnarray}
\begin{aligned}
\int_{\Gamma_1^-}W^{\epsilon_n}V^{\epsilon_n}\be^1\cdot\bn\,ds=&-V^{\epsilon_n}\frac{\mathcal{G}_c}{2}+O(\delta).
\end{aligned}
\label{formulaplusplusfin}
\end{eqnarray}
and \eqref{nonlocSurfaceterm} follows.

To conclude the proof of Proposition \ref{limitspowerflow} we show 
\begin{equation}\label{nonloctolocrate103}
\begin{aligned}
\lim_{{\epsilon_n}\rightarrow 0}E^{\epsilon_n}(\Gamma^{\epsilon_n}_\delta(t))&=-\int_{\Gamma_\delta}\mathbb{C}\mathcal{E}\bu^0\bn\cdot\dot\bu^0\,ds.
 \end{aligned}
\end{equation}
Here we have
\begin{equation}\label{divthm4}
\begin{aligned}
&E^{\epsilon_n}(\Gamma^{\epsilon_n}_\delta(t))\\
&=\epsilon_n^2\int_{P^{\epsilon_n}_\delta(t)}\int_{\mathcal{H}_1(0)}D_{-\be}^{{\epsilon_n}|\xi|}\left[ \epsilon_n|\xi|\partial_S\mathcal{W}^{\epsilon_n}(D_{\be}^{{\epsilon_n}|\xi|}\bu^{\epsilon_n}\cdot\be)\dot\bu^{\epsilon_n}\cdot\be\right]\,d\xi d\bx\\
&=\epsilon_n^2\int_{\mathcal{H}_{1}(0)}\int_{(P^{\epsilon_n}_\delta(t)-\epsilon_n\xi)\setminus P^{\epsilon_n}_\delta(t) }\partial_S\mathcal{W}^{\epsilon_n}(D_{\be}^{{\epsilon_n}|\xi|}\bu^{\epsilon_n}\cdot\be)(\dot\bu^{\epsilon_n}(\bx)+\dot\bu^{\epsilon_n}(\bx+\epsilon_n\xi))\cdot\be\,d\bx d\xi.
\end{aligned}
\end{equation}
Since $\bu^{\epsilon_n}(t)$, $\dot\bu^{\epsilon_n}(t)$, $S(\by,\bx,\bu^{\epsilon_n})$  converge uniformly to $\bu^0(t)$,  $\dot\bu^0$(t), and $\mathcal{E}\bu^0\be\cdot\be$ on $U_\delta(t)$ for $t\in[0,T]$ and $\Gamma^{\epsilon_n}_\delta(t)$ and $\Gamma_\delta$ lie inside  $U_\delta(t)$ we see that
\begin{equation}\label{nonloctolocrate104}
\begin{aligned}
\lim_{{\epsilon_n}\rightarrow 0}E^{\epsilon_n}(\Gamma^{\epsilon_n}_\delta(t))&=\lim_{{\epsilon_n}\rightarrow 0}E^{\epsilon_n}(\Gamma_\delta(t)),
 \end{aligned}
\end{equation}
so we calculate $\lim_{{\epsilon_n}\rightarrow 0}E^{\epsilon_n}(\Gamma_\delta(t))$. Here the contour $\Gamma_\delta(t)$ is the boundary of $P_\delta(t)$ and we write the righthand side of \eqref{divthm4} evaluated now on $P_\delta(t)$ and
\begin{equation}\label{Eongamma}
\begin{aligned}
&E^{\epsilon_n}(\Gamma_\delta(t))\\
&=\epsilon_n^2\int_{\mathcal{H}_{1}(0)}\int_{(P_\delta(t)-\epsilon_n\xi)\setminus P_\delta(t) }\partial_S\mathcal{W}^{\epsilon_n}(D_{\be}^{{\epsilon_n}|\xi|}\bu^{\epsilon_n}\cdot\be)(\dot \bu^{\epsilon_n}(\bx)+\dot\bu^{\epsilon_n}(\bx+\epsilon_n\xi))\cdot\be\,d\bx d\xi\\
&+O(\epsilon_n).
\end{aligned}
\end{equation}

\begin{figure} 
\centering
\begin{tikzpicture}[xscale=0.60,yscale=0.60]

\draw [thick] (-2,-3) rectangle (2,3);

\draw [thick] (-1,-1) rectangle (1,1);


\node [above] at (0,1) {$\Gamma_2$};

\node [right] at (1,0.0) {$\Gamma_3$};

\node [left] at (-1,0) {$\Gamma_1$};

\node [below] at (0,-1) {$\Gamma_4$};

\end{tikzpicture} 
\caption{{ Contour $\Gamma_\delta$ split into four sides.}}
 \label{sides}
\end{figure}
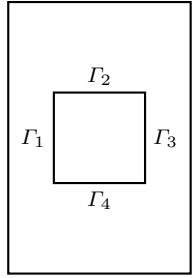

Integration in the $\xi$ variable is over the unit disc centered at the origin $\mathcal{H}_1(0)$. We split the unit disk into its for quadrants $Q_i$, $i=1,\dots,4$. The boundary $\Gamma_\delta$ is the union of its four sides $\Gamma_j$, $j=1,\ldots,4$. Here the left and right sides are $\Gamma_1$ and $\Gamma_3$ respectively and the top and bottom sides are $\Gamma_2$ and $\Gamma_4$ respectively, see Figure \ref{sides}. We choose $\bn$ to be the outward pointing normal vector to $P_\delta$, $\bt$ is the tangent vector to the boundary $\Gamma_\delta$ and points in the clockwise direction, and $\be=\xi/|\xi|$. For $\xi$ in $Q_1$ the set of points $\bx \in (P_\delta(t)-\epsilon_n\xi)\setminus P_\delta(t) $ is parameterized as $\bx=\bt x +\bn(\epsilon_n|\xi|\be\cdot\bn)r$. Here $x$ lies on $\Gamma_1\cup\Gamma_4$ and $0<r<1$ and the area element is $-(\epsilon_n|\xi|\be\cdot\bn)dxdr$. For $\xi$ in $Q_2$ the set of points $\bx \in (P_\delta(t)-\epsilon_n\xi)\setminus P_\delta(t) $ is again parameterized as $\bx=\bt x +\bn(\epsilon_n|\xi|\be\cdot\bn)r$ where $x$ lies on $\Gamma_3\cup\Gamma_4$ and $0<r<1$ and the area element is given by the same formula. For $\xi$ in $Q_3$ we have the same formula for the area element and parameterization and $x$ lies on $\Gamma_3\cup\Gamma_4$  with $0<r<1$. Finally for $\xi$ in $Q_4$ we have again the same formula for the area element and parameterization and $x$ lies on $\Gamma_1\cup\Gamma_2$  with $0<r<1$. This parameterization and a change in order of integration delivers the formula for $E^{\epsilon_n}(\Gamma_\delta(t))$ given by
\begin{equation}\label{Eongammaparam}
\begin{aligned}
&E^{\epsilon_n}(\Gamma_\delta(t))\\
&=-\int_{\Gamma_1}\int_0^1\int_{\mathcal{H}_{1}(0)\cap(Q_1\cup Q_4)}\epsilon_n^3|\xi|\partial_S\mathcal{W}^{\epsilon_n}(D_{\be}^{{\epsilon_n}|\xi|}\bu^{\epsilon_n}\cdot\be)(\dot\bu^{\epsilon_n}(\bx)+\dot\bu^{\epsilon_n}(\bx+\epsilon_n\xi))\cdot\be\bn\cdot\be\,d\xi\,dr\,dx\\
&-\int_{\Gamma_2}\int_0^1\int_{\mathcal{H}_{1}(0)\cap(Q_3\cup Q_4)}\epsilon_n^3|\xi|\partial_S\mathcal{W}^{\epsilon_n}(D_{\be}^{\epsilon_n|\xi|}\bu^{\epsilon_n}\cdot\be)(\dot\bu^{\epsilon_n}(\bx)+\dot\bu^{\epsilon_n}(\bx+\epsilon_n\xi))\cdot\be\bn\cdot\be\,d\xi\,dr\,dx\\
&-\int_{\Gamma_3}\int_0^1\int_{\mathcal{H}_{1}(0)\cap(Q_2\cup Q_3)}\epsilon_n^3|\xi|\partial_S\mathcal{W}^{\epsilon_n}(D_{\be}^{{\epsilon_n}|\xi|}\bu^{\epsilon_n}\cdot\be)(\dot\bu^{\epsilon_n}(\bx)+\dot\bu^{\epsilon_n}(\bx+\epsilon_n\xi))\cdot\be\bn\cdot\be\,d\xi\,dr\,dx\\
&-\int_{\Gamma_4}\int_0^1\int_{\mathcal{H}_{1}(0)\cap(Q_1\cup Q_2)}\epsilon_n^3|\xi|\partial_S\mathcal{W}^{\epsilon_n}(D_{\be}^{{\epsilon_n}|\xi|}\bu^{\epsilon_n}\cdot\be)(\dot\bu^{\epsilon_n}(\bx)+\dot\bu^{\epsilon_n}(\bx+\epsilon_n\xi))\cdot\be\bn\cdot\be\,d\xi\,dr\,dx\\
&+O(\epsilon_n).
\end{aligned}
\end{equation}
Applying the convergence hypothesis of Proposition \ref{limitspowerflow} for $\epsilon_n\rightarrow 0$ one can show using Taylor series that each integrand
converges uniformly to
\begin{equation}\label{limitc}
\begin{aligned}
\frac{4|\xi|}{\omega_2}J(|\xi|)h'(0)\mathcal{E}\bu^0\be\cdot\be(\dot\bu^0\cdot\be)(\bn\cdot\be)
\end{aligned}
\end{equation}
so 
\begin{equation}\label{Eongammaparamlim}
\begin{aligned}
&\lim_{\epsilon_n\rightarrow0}E^{\epsilon_n}(\Gamma_\delta(t))\\
&=-\frac{1}{\omega_2}\int_{\Gamma_1}\int_0^1\int_{\mathcal{H}_{1}(0)\cap(Q_1\cup Q_4)}  4|\xi|J(|\xi|)h'(0)\mathcal{E}\bu^0\be\cdot\be(\dot\bu^0\cdot\be)(\bn\cdot\be)   \,d\xi\,dr\,dx\\
&-\frac{1}{\omega_2}\int_{\Gamma_2}\int_0^1\int_{\mathcal{H}_{1}(0)\cap(Q_3\cup Q_4)}  4|\xi|J(|\xi|)h'(0)\mathcal{E}\bu^0\be\cdot\be(\dot\bu^0\cdot\be)(\bn\cdot\be) \,d\xi\,dr\,dx\\
&-\frac{1}{\omega_2}\int_{\Gamma_3}\int_0^1\int_{\mathcal{H}_{1}(0)\cap(Q_2\cup Q_3)} 4|\xi|J(|\xi|)h'(0)\mathcal{E}\bu^0\be\cdot\be(\dot\bu^0\cdot\be)(\bn\cdot\be)\,d\xi\,dr\,dx\\
&-\frac{1}{\omega_2}\int_{\Gamma_4}\int_0^1\int_{\mathcal{H}_{1}(0)\cap(Q_1\cup Q_2)} 4|\xi|J(|\xi|)h'(0)\mathcal{E}\bu^0\be\cdot\be(\dot\bu^0\cdot\be)(\bn\cdot\be) \,d\xi\,dr\,dx.
\end{aligned}
\end{equation}
Noting that the integrand has radial symmetry in the $\xi$ variable and \eqref{LEFMequality} (see the calculation below Lemma 6.6 of \cite{CMPer-Lipton}) one obtains
 \begin{equation}\label{Eongammaparamlimexplicit}
\begin{aligned}
\lim_{\epsilon_n\rightarrow0}E^{\epsilon_n}(\Gamma_\delta(t))=&-\sum_{i=1}^4\frac{1}{2}\int_{\Gamma_i}2\mathbb{C}\mathcal{E}\bu^0\bn\cdot\dot\bu^0\,dx,
\end{aligned}
\end{equation}
and \eqref{nonloctolocrate103} follows.

\section{Conclusions}
\label{s:conclusions}
In this paper we have shown that the nonlocal cohesive model has solutions that converge in the limit of vanishing non-locality to classic plane elastodynamics with a running crack. The normal traction on the  crack lips  is zero and the energy release rate given by the generalized Irwin relationship (\cite{Freund}, equation (5.39)). The kinetic relation for crack tip motion corresponds to a zero change in internal energy inside domains containing the crack tip and is the classic one given by \eqref{kinetic}.
The  power balance given by \eqref{energyflowbalancedel}, \eqref{energyflowbalance} is not postulated but instead recovered directly from \eqref{energy based model2} by taking the $\epsilon=0$ limit in the nonlocal power balance \eqref{powerbfornonloc}.   In this way one sees that the generalized Irwin relationship is a consequence of the nonlocal cohesive dynamics in the $\epsilon_n=0$ limit.  The recovery is possible since the nonlocal model is well defined over ``the process zone'' around the crack centerline tip. This shows that the double well potential provides a phenomenological description of the process zone at mesoscopic length scales.  

In this paper we have illustrated the ideas using the simplest double well energy for a bond based perydynamic formulation.  We are free to take a more sophisticated energy like those motivated by the Lennard Jones potential. Doing so will deliver a nonlocal model that preserves non-interpenetration of material points for all types of loadings. We can then pass to the small horizon limit in such a model to recover a sharp fracture model with crack lips that do not interpenetrate. More generally we may consider state based peridynamic models and perform similar analyses. These are projects for the future but all are theoretically accessible.

\newcommand{\noopsort}[1]{}

\end{document}